%% file: 2005-17.tex
\documentclass{gtart_h}  

\input gtspec

\lognumber{510}
\received{31 October 2004}
\volumenumber{9}\papernumber{17}\volumeyear{2005}
\pagenumbers{699}{755}   
\revised{18 April 2005}
\published{28 April 2005}
\accepted{6 March 2005}
\proposed{Robion Kirby}
\seconded{Ronald Fintushel,  Ronald Stern}

\usepackage{amsmath, amssymb}
\usepackage[curve]{xypic}

\newcommand{\comment}[1]{}
\newbox\mybox
\def\overtag#1#2#3{\setbox\mybox\hbox{$#1$}\hbox to
  0pt{\vbox to 0pt{\vglue-#3\vglue-\ht\mybox\hbox to \wd\mybox
      {\hss$\ss#2$\hss}\vss}\hss}\box\mybox}
\def\undertag#1#2#3{\setbox\mybox\hbox{$#1$}\hbox to 0pt{\vbox to
    0pt{\vglue#3\vglue\ht\mybox\hbox to \wd\mybox
      {\hss$\ss#2$\hss}\vss}\hss}\box\mybox}
\def\lefttag#1#2#3{\hbox to 0pt{\vbox to 0pt{\vss\hbox to
      0pt{\hss$\ss#2$\hskip#3}\vss}}#1}
\def\righttag#1#2#3{\hbox to 0pt{\vbox to 0pt{\vss\hbox to
      0pt{\hskip#3$\ss#2$\hss}\vss}}#1}
\let\ss\scriptstyle

\def\Dot{\lower.2pc\hbox to 2.5pt{\hss$\bullet$\hss}}
\def\Circ{\lower.2pc\hbox to 2.5pt{\hss$\circ$\hss}}
\def\Vdots{\raise5pt\hbox{$\vdots$}}
\def\splicediag#1#2{\xymatrix@R=#1pt@C=#2pt@M=0pt@W=0pt@H=0pt}

\renewcommand\frame[2][3pt]{\hbox{$\vcenter{\hbox{\vrule\vbox
{\hrule\kern#1\hbox{\kern#1$#2$\kern#1}\kern#1\hrule}\vrule}}$}}

\def\minimal{quasi-minimal}\def\Minimal{Quasi-minimal}
\def\z{B}
\def\disc{\det}
\newcommand\lineto{\ar@{-}}
\newcommand\dashto{\ar@{--}}
\newcommand\dotto{\ar@{.}}
\newcommand{\C}{{\mathbb C}}
\newcommand{\Z}{{\mathbb Z}}
\newcommand{\N}{{\mathbb N}}
\newcommand{\Q}{{\mathbb Q}}
\newcommand{\E}{{\mathbb E}}
\newcommand{\resgraph}{\Gamma}
\newcommand{\lcm}{\operatorname{lcm}}
\newcommand{\Hom}{\operatorname{Hom}}
\newtheorem{theorem}{Theorem}[section]
\newtheorem*{theorem*}{Theorem}
\newtheorem{lemma}[theorem]{Lemma}
\newtheorem{proposition}[theorem]{Proposition}
\newtheorem{corollary}[theorem]{Corollary}

\newtheorem{conjecture}[theorem]{Conjecture}
\newtheorem*{condition3.3}{`Condition 3.3'}
\newtheorem*{condition3.4}{`Condition 3.4'}
\theoremstyle{definition}
\newtheorem{definition}[theorem]{Definition}
\newtheorem{example}[theorem]{Example}

\newtheorem{notation}[theorem]{Notation}

\begin{document}

\title{Complete intersection singularities of splice type\\as
  universal abelian covers}
\shorttitle{Complete intersection singularities of splice type}
\author{Walter D Neumann\\Jonathan Wahl}
\address{Department of Mathematics, Barnard College, Columbia
  University\\New York, NY 10027, USA}
\address{Department of Mathematics, The University of
North Carolina\\Chapel Hill, NC 27599-3250, USA}
\asciiaddress{Department of Mathematics, Barnard College, Columbia
  University\\New York, NY 10027, USA\\and\\Department of 
Mathematics, The University of
North Carolina\\Chapel Hill, NC 27599-3250, USA}
\gtemail{\mailto{neumann@math.columbia.edu}{\rm\qua 
and\qua}\mailto{jmwahl@email.unc.edu}}
\asciiemail{neumann@math.columbia.edu, jmwahl@email.unc.edu}

\asciikeywords{Surface singularity, Gorenstein singularity,
rational homology sphere, complete intersection singularity,
abelian cover}

\begin{abstract}
  It has long been known that every quasi-homogeneous normal complex
  surface singularity with $\Q$--homology sphere link has universal
  abelian cover a Brieskorn complete intersection singularity. We
  describe a broad generalization: First, one has a class of complete
  intersection normal complex surface singularities called ``splice
  type singularities,'' which generalize Brieskorn complete intersections.
  Second, these arise as universal abelian covers of a class of normal
  surface singularities with $\Q$--homology sphere links, called
  ``splice-quotient singularities.''  According to the Main Theorem,
  splice-quotients realize a large
  portion of the possible topologies of singularities with
  $\Q$--homology sphere links.  As quotients of complete intersections,
  they are necessarily $\Q$--Gorenstein, and many $\Q$--Gorenstein
  singularities with $\Q$--homology sphere links are of this type. We
  conjecture that
  rational singularities and minimally elliptic singularities with
  $\Q$--homology sphere links are splice-quotients. A recent preprint
  of T Okuma presents confirmation of this conjecture.
\end{abstract}
\asciiabstract{%
It has long been known that every quasi-homogeneous normal complex
  surface singularity with Q-homology sphere link has universal
  abelian cover a Brieskorn complete intersection singularity. We
  describe a broad generalization: First, one has a class of complete
  intersection normal complex surface singularities called "splice
  type singularities", which generalize Brieskorn complete
  intersections.  Second, these arise as universal abelian covers of a
  class of normal surface singularities with Q-homology sphere
  links, called "splice-quotient singularities".  According to the
  Main Theorem, splice-quotients realize a large portion of the
  possible topologies of singularities with Q-homology sphere
  links.  As quotients of complete intersections, they are necessarily
  Q-Gorenstein, and many Q-Gorenstein singularities with
  Q-homology sphere links are of this type. We conjecture that
  rational singularities and minimally elliptic singularities with
  Q-homology sphere links are splice-quotients. A recent preprint
  of T Okuma presents confirmation of this conjecture.}

\maketitle

\section{Introduction}
The possible topologies for a normal singularity of a complex surface
are classified (eg, \cite{neumann81}), but it is very
rare that much is known about possible analytic types for given
topology.  Locally, the topology is the cone on an oriented 3-manifold
$\Sigma$, called the \emph{link} of the singularity.  Via the
configuration of exceptional curves on a good
resolution $(\tilde Y ,E)\rightarrow (Y,o)$, one can construct
$\Sigma$ via plumbing according to the negative-definite dual
resolution graph $\Gamma$.  In this paper we will restrict to the case
$\Sigma$ is a \emph {rational homology sphere}, or $\Q$HS, ie,
one for which
$H_{1}(\Sigma;\Z)$
is finite; equivalently, the exceptional configuration is a tree of
smooth rational curves.
%This enormous class includes all rational and
%nearly all minimally elliptic singularities, the hypersurfaces
%$\{x_{p}+y^{q}+z^{r}=0\}$ if $(p,qr)=1$, etc.

The universal abelian covering
$\tilde \Sigma \rightarrow \Sigma$
is finite, and can be realized by a finite map of germs $(X,o)\rightarrow
(Y,o)$; the covering (or discriminant) group
$H_{1}(\Sigma)$ is easily computed from the dual graph
$\Gamma$.  Given $\Sigma$, or (equivalently, by \cite{neumann81})
a graph $\Gamma$, our goal is to
construct an explicit singularity $(X,o)$ whose link is $\tilde \Sigma$, and an
action of the discriminant group which is free off $o$, so that the
quotient $(Y,o)$ has graph $\Gamma$.  We will achieve this under
certain conditions on $\Gamma$ (Theorem \ref{th:main}).

Suppose first that $\Sigma$ is Seifert fiberable.  Then it has been
known for some time that the universal abelian cover of $\Sigma$ is
diffeomorphic to the link of a Brieskorn complete intersection
singularity\footnote{In this paper diffeomorphisms are always assumed
  to preserve orientation; and, since we are interested in
  singularities, complete intersections are local complete
  intersections in the usual sense -- eg, \cite{hartshorne}, p.\  185.}
(\cite{neumann81, neumann83}). Thus a possible analytic type is as an
abelian quotient of a Brieskorn complete intersection.  From another
point of view, consider a quasi-homogeneous $(Y,o)$ with $\Q$HS link;
the resolution diagram is star-shaped, and from it one can read off
easily the data needed to write down a Brieskorn complete intersection
and a diagonal action of the discriminant group.  A look at the
Seifert data shows that the quotient has the same topology as $(Y,o)$.
As a bonus, one can even arrange to recover the analytic type of $Y$,
because one knows the exact ingredients needed to make a
quasi-homogeneous singularity.

At this point, to handle more general links, one might wonder what kinds
of equations could
generalize Brieskorn complete intersections.

Three-manifold theory gives a natural minimal decomposition of
$\Sigma$ along embedded tori into pieces that are Seifert fiberable
(a version of the JSJ-decompos\-ition,  \cite{neumann-swarup}).  When
  this set of tori is empty, $\Sigma$ is the only piece in the
  decompos\-ition---the aforementioned Seifert case.  More generally,
  associated to the JSJ decomposition of $\Sigma$ is a certain
  weighted tree called a \emph{splice diagram}. This is a tree with no
  valence two vertices, and for each node (vertex of valence $\ge 3$)
  it has a positive integer weight associated with each incident edge.
  The pieces in the JSJ decomposition of $\Sigma$ are in one-one
  correspondence with the nodes of the splice diagram. The splice
  diagram does not necessarily determine $\Sigma$ but it does
  determine its universal abelian cover.

As indicated, the link $\Sigma$ determines the topology of the minimal good
resolution of the singularity, and the splice diagram $\Delta$ can
easily be computed from the resolution dual graph $\Gamma$
(see, eg, Section \ref{sec:splicing}). $\Delta$
has the same general shape as $\Gamma$, but degree two vertices are
suppressed. Resolution graphs satisfy a negative definiteness
condition which translates into a ``positive edge
determinants'' condition for the splice diagrams of
singularities.

Under a certain natural condition on $\Delta$, called the ``semigroup
condition,'' we associate a collection of $t-2$ equations in $t$
variables, where $t$ is the number of leaves of $\Delta$. There
is some choice allowed in these ``splice diagram equations.''
They generalize Brieskorn complete intersections as
follows. A Brieskorn complete intersection, corresponding to a splice
diagram with a single node of valence $t$, is defined by a system of
$t-2$ weighted homogeneous equations. For a splice diagram with more than
one node, one associates to each node a collection of $\delta-2$ equations
($\delta$ the valence of the node) which are weighted homogeneous with
respect to a system of weights associated to the node.  This gives a
total of $t-2$ equations. (We also allow higher weight perturbations
of these equations.)

 We also formulate the ``congruence
conditions,'' which depend on $\Gamma$, which guarantee that the
discriminant group of the resolution acts on a set of splice diagram
equations for $\Delta$.  Our main results (Theorems \ref{th:splice is
  CI} and Theorem \ref{th:main}) can be summarized:
\begin{theorem*}
The splice diagram equations associated to a splice diagram $\Delta$ with
semigroup condition always describe a normal complete
intersection singularity.

If splice diagram equations have been chosen equivariantly with
respect to the action of the discriminant group for $\Gamma$, then
the action is free away from the singular point, the quotient is a
singularity with resolution graph $\Gamma$ (and hence with link
$\Sigma$), and the covering is the universal abelian cover.
\end{theorem*}
\noindent Thus for a large
family of topologies --- those that satisfy the semigroup and
congruence conditions --- we find explicit (and attractive) analytic
descriptions of singularities with the given topology.
Put another way, in such cases we can write down
explicit equations for singularities with given topology (of course,
modulo writing down invariants for the group action).

As finite quotients of complete intersections, these ``splice-quotient
singularities'' are necessarily $\Q$--Gorenstein. Although splice
diagram equations depend on choice of certain ``admissible
monomials,'' the family of analytic types for the resulting
splice-quotient singularities is independent of these choices (Theorem
\ref{th:change monomial})\comment{addition}.  We had earlier (rashly)
conjectured that every $\Q$--Gorenstein singularity with $\Q$HS link
should be a splice-quotient, and though this appears to be true
surprisingly often (eg, quasihomogeneous singularities
\cite{neumann83} and quotient-cusps \cite{neumann-wahl03}),
counterexamples are now known \cite{nemethi et al}. In fact, weakly
equisingular deformations (in the sense of weak simultaneous
resolution) of a splice-quotient need not be of that type.
Nevertheless, we conjectured in the original version of this paper
that rational singularities and minimally elliptic singularities with
$\Q$HS link should be splice quotients (Conjecture \ref{rational}).
The recent preprint \cite{okuma} of T. Okuma now offers a proof of
this Conjecture (see Section \ref{sec:okuma}).

In the important case that $\Sigma$ is a $\Z$--homology sphere (these
are classified by splice diagrams with pairwise prime weights at each
node)\comment{Changed: being explicit was hardly longer}, no abelian
quotient is needed; so, the semigroup condition then implies one can
write down directly complete intersection splice equations with given
topology.  We know of no complete intersection with $\Z$--homology
sphere link which does not satisfy the semi-group conditions, or is
not of splice type.  (See \cite{neumann-wahl-zsphere}.)

The leaves of the splice diagram correspond to knots in $\Sigma$, and
we show that they are cut out (in the universal abelian cover) by
setting the corresponding variable equal to zero. In
\cite{neumann-wahl-zsphere} we show in the $\Z$--homology sphere case
that the existence of functions cutting out these knots is equivalent
to the singularity being of splice type, and we conjecture that this
holds more generally%\footnote{This is also discussed further in
%  Section \ref{sec:okuma} on Okuma's preprint.}
. This point of view is
useful beyond the question of existence of singularities with given
topologies---it applies as well to analytic realization of germs of
curves in complex surfaces (especially with $\Z$--homology sphere
links).  For instance, embedded resolution of a plane curve
singularity gives rise to a (non-minimal) $\Gamma$ and $\Delta$, and
one writes down an explicit equation of the curve by setting a
variable equal to $0$, as in \cite{neumann-wahl-zsphere}, Section 5.
More generally, consider the following illustration. Let $X$ be the
Brieskorn variety $x^2+y^3+z^{13}=0$, and let $K$ be the knot in its
link cut out by $z=0$ (this is the degree $13$ fiber of the Seifert
fibration of $\Sigma$).  Form $K(p,q)$, the $(p,q)$--cable on $K$, a
new knot on $\Sigma$.  Then the positive edge determinant condition
says that $K(p,q)$ is the link of a complex curve through the origin
in $X$ if and only if $13q>6p$, and the semigroup condition says that
this curve can be cut out by a single equation $f(x,y,z)=0$ if (and in
this case only if) $q\ge2$.  We will return to this theme elsewhere.

\smallskip Let us explain the steps needed to get to the main result.
First, given a splice diagram $\Delta$ satisfying the semigroup
conditions, we write down an explicit set of equations, and our first
goal is to show (Theorem \ref{th:splice is CI}) that these splice
diagram equations define an isolated complete intersection
singularity.  Every node of $\Delta$ defines a weight filtration, and
one needs to prove that each associated graded is a reduced complete
intersection, defined by the leading forms of the given equations.
This step (Section \ref{sec:assoc graded}) involves understanding
curves defined by analogs of splice diagram equations and some
detailed combinatorics involving the diagram weights.  Then, to show
the singularity is isolated\comment{Added phrase}
(Section \ref{sec:blow-up}), one does a
weighted blow-up at an ``end-node'' of $\Delta$, and examines
singularities along the exceptional fiber.  The key is to show one now
has equations for a ``smaller'' $\tilde \Delta$, which by induction
has an isolated singularity; a difficulty is that the new weights are
related to the old ones in a rather complicated way.

We next consider the resolution diagram $\Gamma$. The discriminant
group $D(\Gamma)$ is computed, and shown to act naturally (and without
pseudoreflections) on $\C ^{t}$, the space on which the splice diagram
equations for $\Delta$ are defined (Section \ref{sec:discriminant}).
To proceed, one needs to be able to choose splice diagram equations on
which $D(\Gamma)$ acts equivariantly; this is the ``congruence
condition'' which we need.  In particular, if the semigroup and
congruence conditions are satisfied for $\Gamma$, then we have an
action of the discriminant group on the splice diagram equations. Our
main theorem, Theorem \ref{th:main}, asserts that the quotient map is
the universal abelian cover of a normal singularity whose resolution
dual graph is $\Gamma$.

Finally, we have to prove Theorem \ref{th:main}, which we do by
induction on the number of nodes.  A key point is to explicitly lift
the generators of $D(\Gamma)$ to a weighted blow-up of the singularity
at an ``end-node''.  We then have to identify the exceptional fiber
and singularities that arise after factoring by the lifted group.  We
show that at the ``worst'' singular points one has (after analytic
change of coordinates) a splice-quotient for a subgraph
$\tilde \Gamma$ of $\Gamma$.  Again, this involves some complicated
numerics, some of which are proved in the Appendix. The one-node case
is done in Section \ref{sec:one node} and the inductive step in
Section \ref{sec:induct}.

We also show, as part of the main theorem, that for any node of the
splice diagram, the grading on the splice-quotient induced by the
weight filtration on $\C^t$ is, up to a multiple, just order of
vanishing on the corresponding curve of the
resolution\comment{Addition}. Using this we can deduce (Theorem
\ref{th:change monomial}) that the concept of splice-quotient
singularity is a canonical concept (independent of choices).

In Section \ref{sec:twonode}, we take a two-node minimal resolution graph, and
write down explicitly the semigroup and congruence conditions.
\iffalse We
show that in the rational and $\Q$HS minimally elliptic cases,
these conditions are automatic, as is consistent with the much stronger
Conjecture \ref{rational}.\fi

The first Appendix (Section \ref{sec:splicing}) proves some results about
resolution and splice diagrams that are needed in the paper, as well
as a topological description of splice diagrams.

A second Appendix, added April 2005, discusses Okuma's recent preprint
on the conjecture that rational and $\Q$HS-link minimally
elliptic singularities are always splice quotients.

{\bf Acknowledgements}\qua The first author's research is supported
under NSF grant DMS-0083097 and the second author's under NSA
grant MDA904-02-1-0068.

\section{Semigroup conditions and splice equations}\label{sec:basics}

We first recall the concept of ``splice diagram,'' a certain kind of
weighted tree.  Given a finite tree, the \emph{valency} of a vertex is
the number of incident edges. A \emph{node} is a vertex of valency
$\ge 3$ and a \emph{leaf} is a vertex of valency $1$.

A \emph{splice diagram} $\Delta$ is a finite tree with no valence $2$
vertices, decorated with integer weights as follows: for each node $v$
and edge $e$ incident at $v$ an integer weight $d_{ve}$ is given (see
Section \ref{sec:splicing} for examples and more
detail). Thus an
edge joining two nodes has weights associated to each end, while an
edge from a node to a leaf has just one weight at the node end.  The
\emph{edge determinant} of an edge joining two nodes is the product of
the two weights on the edge minus the product of the weights adjacent
to the edge.  Splice diagrams that arise in the
study of links of complex singularities always satisfy the following
conditions:
\begin{itemize}
\item All weights are positive.
\item All edge determinants are positive.
\item One has the ``ideal condition'' on weights (see below).
\end{itemize}
To explain the third of these, we need more notation.

\begin{notation}\label{notation2}
  For a node $v$ and an edge $e$ at $v$, let $d_{ve}$ be the weight on
  $e$ at $v$, and $d_{v}$ the product of all edge-weights $d_{ve}$
  at $v$.
   For any pair of distinct vertices $v$ and $w$, let $\ell_{vw}$ (the
  \emph{linking number}) be the product of all the weights adjacent
  to, but not on, the shortest path from $v$ to $w$ in $\Delta$.
  Define $\ell'_{vw}$ similarly, except one
  excludes weights
  around $v$ and $w$. (Thus $\ell'_{vw}=1$ if $v$ and $w$ are
  adjacent, and $\ell_{vw}=\ell'_{vw}$ if $v$ and $w$ are both
  leaves.)
  Finally, let $\Delta_{ve}$ be the subgraph of $\Delta$ cut off from $v$ by
  $e$ (ie, on the ``$e$--side of $v$'').
\end{notation}

\begin{definition}[Ideal Condition]\label{def:ideal}
  For each node $v$ and adjacent edge $e$ of $\Delta$, the edge-weight
  $d_{ve}$ is divisible by the GCD of all $\ell'_{vw}$ with $w$ a leaf
  of $\Delta$ in $\Delta_{ve}$; in other words, $d_{ve}$ is in the
  following ideal of $\Z$:
$$d_{ve}\in\left(\ell'_{vw}: w \text{ a leaf of }\Delta\text{ in }
  \Delta_{ve}\right)\,.$$
\end{definition}

We are interested in splice diagrams that satisfy the stronger  condition:

\begin{definition}[Semigroup Condition]\label{def:sg}
  The \emph{semigroup condition} says that for each node $v$ and
  adjacent edge $e$ of $\Delta$, the edge-weight $d_{ve}$ is in the
  following sub-semigroup of $\N$:
$$d_{ve}\in\N \langle \ell'_{vw}:w \text{ a leaf of $\Delta$ in }
  \Delta_{ve} \rangle \,.$$
\end{definition}

For each edge $e$ at $v$ the semigroup condition lets us write
\begin{subequations}\label{eq:sg}
\begin{equation}\label{eq:sg1}
d_{ve}=\sum_{w\text{ a leaf in
  }\Delta_{ve}}\alpha_{vw}\ell'_{vw}\,,\quad\text{with }
\alpha_{vw}\in\N \cup\{0\}.
\end{equation}
It is easy to see that this is equivalent to
\begin{equation}\label{eq:sg2}
d_{v}=\sum_{w}\alpha_{vw}\ell_{vw}.
\end{equation}
\end{subequations}
Assume from now on that $\Delta$ satisfies the semigroup condition.
To each leaf $w$ we associate a variable $z_{w}$.

\begin{definition}[$v$--weighting; admissible
  monomials]\label{admissible}
  Fix a node $v$. Then the \emph{$v$--weighting}, or
  \emph{$v$--filtration}, of the polynomial ring in the $z_{w}$'s is
  defined by assigning weight $\ell_{vw}$ to $z_{w}$.

An \emph{admissible monomial} (associated to
the edge $e$ at the node $v$) is a monomial $M_{ve}=\prod_w
z_w^{\alpha_{vw}}$, the product over leaves $w$ in $\Delta_{ve}$, with
exponents satisfying the above equations (\ref{eq:sg}).
In particular, each admissible monomial
$M_{ve}$ is $v$--weighted homogeneous, of total $v$--weight $d_{v}$.
\end{definition}

\begin{definition}[Splice diagram equations]  Let $\Delta$ be a
  splice diagram with $t$ leaves satisfying the semigroup condition.
  To each leaf $w$ associate a variable $z_{w}$; for each node $v$ and
  adjacent edge $e$, choose an admissible monomial $M_{ve}$.  Then
  \emph{splice diagram equations for $\Delta$} consist of a collection
  of equations of the form
  $$\sum_{e}a_{vie}M_{ve} +H_{vi}=0,\quad v\text{ a node},~
  i=1,\dots,\delta_{v}-2\,,$$
  where
\begin{itemize}
\item for every $v$, all maximal minors of the
  $((\delta_v-2)\times\delta_v)$--matrix $(a_{vie})$ have full rank
\item $H_{vi}$ is a convergent power series in the $z_{w}$'s all of
  whose monomials have $v$--weight $>d_{v}$.
\end{itemize}

It is easy to see one has exactly $t-2$ equations in the $t$ variables.
The corresponding
subscheme $X(\Delta)\subset \C^{t}$ is a
\emph{splice diagram surface singularity}.
\end{definition}

\begin{theorem}\label{th:splice is CI} Let $X=X(\Delta)\subset \C^{t}$ be a splice diagram
surface singularity.  Then:
\begin{enumerate}
\item $X$ is a two-dimensional complete
intersection with an isolated singularity at the origin.
\item For any node $v$, the corresponding weight filtration has
associated graded ring a reduced complete intersection, defined by
the $v$--leading forms of the splice equations.
\end{enumerate}
\end{theorem}

The theorem will be proved partly by induction on the number of nodes
of $\Delta$.  Once we know the singularities are isolated, one can
recover all the analytic types by restricting to polynomials $H_{vi}$
in the definition. We allow splice diagrams in which an edge-weight
leading to a leaf may be 1, as
the ``minimality'' assumption which
avoids this could be lost in our inductive process.

One could define a more general class of equations by allowing, for
each node $v$ and edge $e$, linear combinations of admissible
monomials, rather than multiples of a fixed one.  It follows from the
theorem (and proof) that for \emph{generic coefficients}, these give
isolated singularities with the same properties.  But, in the
situation of most interest to us here (Theorem \ref{th:main}), where
the monomial $M_{ve}$ also satisfies an equivariance condition, this
is not a generalization, since we will show (Theorem \ref{th:change
  monomial}) that any other allowed monomial $M'_{ve}$ for $v,e$ then
differs from some multiple of $M_{ve}$ by something of higher weight.

Given that the maximal minors of the
coefficient matrix $(a_{vie})$ have full rank, one may apply row
operations to the matrix (which is the same as taking linear
combinations of the corresponding equations), to put the
$(\delta_v-2)\times\delta_v$ coefficient matrix in the form
$$
\begin{pmatrix}
1&0&\dots&0&a_1&b_1  \\
0&1&\dots&0&a_2&b_2\\
\vdots&\vdots&\ddots&\vdots&\vdots&\vdots\\
0&0&\dots&1&a_{\delta_v-2}&b_{\delta_v-2}
\end{pmatrix}
$$
with $a_ib_j-a_jb_i\ne0$ for all $i\ne j$, and all $a_i$ and $b_i$
nonzero.  We will often assume we have done so. In this way, the
defining equations are sums of three monomials, plus higher order
terms. (Sometimes --- for instance the next section --- it will be
more convenient to move the last two columns of the above matrix to
the first two columns.)

\section{Splice diagram curves and associated gradeds}\label{sec:assoc graded}

As in \cite{neumann-wahl-zsphere}, one can define curves using a
modified version of splice diagrams.  Let $(\Delta,w')$ be a splice
diagram with distinguished leaf $w'$. But now, at any node, the edge
weight in the direction of $w'$ is irrelevant and should be omitted or
ignored.  $(\Delta,w')$ \emph{satisfies the semigroup condition} if
for every node $v$ and adjacent edge $e$ pointing away from $w'$, the
edge-weight $d_{ve}$ is in the semigroup generated by all the
$\ell'_{vw}$ with $w$ in $\Delta_{ve}$.  This condition is automatic
if $\Delta$ is a splice diagram that satisfies the semigroup condition
and $(\Delta,w')$ is obtained by ignoring edge weights in the
direction of $w'$.

As before, associate a variable $z_{w}$ to each leaf different from
$w'$; and, for every node $v$ and adjacent edge $e$ pointing away from
$w'$, choose an admissible monomial $M_{ve}$.  One can define a
\emph{splice diagram curve $C=C(\Delta,w')$} via equations
for all nodes $v$:
$$\sum_{e}a_{vie}M_{ve}=0,\quad i=1,\ldots,\delta_{v}-2\,,$$
where for all
$v$ the $(\delta_{v}-2)\times (\delta_{v}-1)$ matrix $(a_{vie})$ is
required to have full rank.  (Note that we now don't allow higher
order terms.)  Enumerating the admissible monomials at
$v$ as $M_{vj}$, $1\leq j \leq \delta_{v}-1$, one can always write the
system of equations as
$$M_{vj}=a_{vj}M_{v1},\quad j=2,\dots,\delta_{v}-1\,,$$
where all
$a_{vj}\neq 0$. Assign to the variable $z_{w}$ the weight
$\ell_{w'w}$; then the equations at $v$ are weighted homogeneous, of
total weight $\ell_{w'v}$.

\begin{theorem}\label{th:curves} Let $C$ be a splice diagram curve as above.
\begin{enumerate}
\item At a point of $C$ for which one coordinate is 0, all coordinates
  are 0.
\item Except at the origin, $C$ is a smooth curve.
\item The number of components of $C$ is the GCD of the weights of the
  variables.
\item $C$ is a reduced complete intersection curve, and every
  irreducible component is isomorphic to a monomial curve.
\item For each leaf $w\neq w'$, let $\alpha_{w}$ be an integer $\geq
  0$, not all equal 0.  Then for $b\neq 0$, the intersection
  $C\cap\{\prod z_{w}^{\alpha_{w}}=b\}$
  is everywhere transverse and
  consists of\/ $\sum \alpha_{w}\ell_{w'w}$ points.
\end{enumerate}
\end{theorem}
\begin{proof}  Denote by $v'$ the node adjacent to the distinguished
  leaf $w'$.

  One can use induction in several ways: for instance, remove $v'$
  from the diagram, then reinsert a distinguished leaf in its place
  for $\delta_{v'}-1$ splice subdiagrams pointing away from $w'$. The
  splice type equations for the subdiagrams are among the equations we
  started with.  Combining with the equations at $v'$, one easily gets
  the first claim.

  For the second assertion, it is easier to use a different induction.
  Let $v$ be an end-node, adjacent to $\sigma:=\delta_v-1$
  leaves, say $w_{1},\dots,w_{\sigma}$, with associated coordinates
  $z_{1},\dots,z_{\sigma}$.  Then some splice equation at $v$ is of
  the form $z_{2}^{m_{2}}-
  az_{1}^{m_{1}}=0$ ($a\neq 0$); this polynomial has $k$
  irreducible factors, where $k$ is the GCD of $m_{1}$ and
  $m_{2}$.  In a neighborhood of the point in question, exactly
  one of these factors vanishes, so we can solve
  $z_{1}=t^{m_{2}/k},
  z_{2}=a't^{m_{1}/k}$.  Plug these into all the
  remaining equations.  One then recognizes splice diagram equations
  associated to a diagram with one less leaf.  When $\sigma \geq 3$,
  replace $w_{1}$ and $w_{2}$ by one leaf, with
  edge-weight $m_{1}m_{2}/k$, and then divide all edge
  weights pointing towards the node $v$ by $h$ (they are divisible by
  $k$, by the semigroup condition).  When $\sigma =2$, drop the two
  leaves, making $v$ into a new leaf, and again divide the same
  edge-weights by $k$.  It is straightforward to check that one has
  equations associated to a new splice diagram.  This process gives
  the inductive step necessary to prove the claim.

  The third assertion uses the same induction.
Number the leaves
  $w_1,\dots,w_\tau$ and abbreviate the weight $\ell_{w'w_i}$ of the
  $i$-th variable by $\ell_i$. Then the weights of the variables in
  the new splice diagram in the induction are
  $\ell:=\ell_1/m_2=\ell_2/m_1$, and $\ell_3/k, \dots,\ell_\tau/k$.  The
  inductive step thus follows from the equation
  $$\gcd(m_1\ell,m_2\ell,\ell_3,\dots,\ell_\tau)=
  k.\gcd(\ell,\ell_3/k,\dots,\ell_\tau/k)\,.$$
  Since $C$ is a curve with one singularity, and is defined by the
  appropriate number of equations, it is a complete intersection
  (necessarily reduced).  Since $C$ is weighted homogeneous, it has a
  $\C^*$--action. Any irreducible component still has a
  $\C^{*}$--action, hence is a monomial curve.  In fact, let $h$ denote
  the GCD of all the weights $\ell_i$ of the variables, and consider
  ``reduced weights'' $\bar\ell_{i}:= \ell_{i}/h$. \comment{Is
  $\bar{\ell_{i}}$ better? $\ell '$ already used. DONE}
  Then for any
  $(p_{1},\dots,p_{\tau})\neq (0,\dots,0)$ on $C$, one sees that
  $(z_{1},\dots,z_\tau)=(p_{1}t^{\bar\ell_{1}},\dots,p_{\tau}t^{\bar\ell_{\tau}})$
  defines an irreducible monomial curve contained in $C$ through the
  point, hence gives an irreducible component.

  As for the last count of solutions,
  the function $f:=\prod z_w^{a_w}$, restricted to the component
  $(p_{1}t^{\bar\ell_{1}},\dots,p_{\tau}t^{\bar\ell_{\tau}})$ of $C$, has
  the form $pt^{\sum a_i\bar\ell_i}$ for some $p\ne0$, so it has exactly
  $\sum a_i\bar\ell_i$ inverse images over any $b\ne0$. Thus, the
  intersection of $C$ with $f(z_1,\dots,z_\tau)=b$ is transversal for
  any $b\ne 0$ and the number of intersection points is $\sum
  a_i\bar\ell_i$ on each of the $h$ components of $C$, so there are
  $h\sum a_i\bar\ell_i=\sum a_i\ell_i$ intersection points in all.
\end{proof}

 From the preceding result one derives the second part of Theorem
\ref{th:splice is CI} concerning $X(\Delta)$.
\begin{proof}[Proof of Theorem \ref{th:splice is CI} (2)]
  Enumerate the chosen admissible monomials at $v$ as $M_{i}$, $1\leq
  i \leq \delta_{v}$; these are weighted homogeneous of degree $d_{v}$
  with respect to the $v$--weighting.  So the leading forms of the
  equations at $v$ may be written
$$M_{i}-a_{i}M_{1}-b_{i}M_{2}=0,\quad  3\leq i
\leq \delta_{v}\,.$$

\begin{lemma}\label{le:weights}
At a node $v'\neq v$, consider an admissible monomial $M_{v'e}$.
\begin{enumerate}  \item If $e$ does not point towards $v$, then the
$v$--weight of $M_{v'e}$ is $\ell_{v'v}$.
\item If $e$ points towards $v$, then the $v$--weight of $M_{v'e}$
is strictly greater than $\ell_{v'v}$.
\item If a monomial has $v'$--weight $>d_{v'}$, then its
$v$--weight is strictly greater than $\ell_{v'v}$.
\end{enumerate}
\end{lemma}
\begin{proof} (Cf.\ Theorem 7.3 in
    \cite{neumann-wahl-zsphere}.)
One checks directly that for any $w$
  $$\frac{\ell_{vw}}{\ell_{v'w}}=\frac{\ell_{v'v}}{d_{v'}}D_{e'_1}\dots
D_{e'_k}$$
where:
\begin{itemize}
\item $e'_1,\dots,e'_k$ are the edges that are on the path from $v'$
  to $v$ but not on the path from $w$ to $v$;
\item   for any edge $e$, $D_e$ is the product of the edge
  weights on $e$ divided by the product of edge weights directly
  adjacent to $e$ (so $D_e>1$ by the edge determinant condition).
\end{itemize}
Thus $\ell_{vw}/\ell_{v'w}$ takes its minimum value (namely
${\ell_{v'v}}/{d_{v'}}$) if and only if $w$ is beyond $v'$ from the
point of view of $v$.  It follows that the admissible monomials at
$v'$ all have the same $v$--weight except for the near monomial for
$v$, which has higher $v$--weight. One has the last assertion of the
lemma as well.
\end{proof}

It follows from the lemma that the $v$--leading form of a splice
diagram equation at $v'$ is given by the combination of the
admissible monomials pointing away from $v$; all other terms
have higher $v$--weight.  We prove that these equations, plus the
ones at $v$, together define a two-dimensional complete intersection,
whence they define the full associated graded ring associated to the
$v$--filtration.

The vertices $\neq v$ of $\Delta$ divide into $\delta_{v}$ groups,
depending upon the $\Delta_{ve}$ to which a $v'$ belongs.
In a
sector consisting of more than one leaf, consider the
$v$--leading forms
 of the splice equations corresponding to vertices in that sector.
Using Lemma \ref{le:weights}, these forms can be seen to provide a complete set
 of splice diagram curve equations for the corresponding set of variables
 (where $v$ is viewed as a ``root'' of the sector).  Therefore, in any sector,
 if one of the $z_{w}$ vanishes at
a point, so do all the other leaf variables (by Theorem
\ref{th:curves}).
In particular,
at a point of vanishing of two variables, one from each
of the first two edge directions, one has the vanishing as well of
$M_{1}$ and $M_{2}$.
This implies that every $M_{i}$ vanishes, whence one variable in every
group vanishes, so every variable vanishes.
Consequently, the locus
of the particular $t-2$ leading forms plus these two variables is
0-dimensional, hence must be a 0-dimensional complete intersection.
Equivalently, those leading forms and two variables form a regular
sequence.  Now recall the following well-known result in
commutative algebra:

\begin{lemma}  Let $f_{1},\ldots ,f_{r}$ be a sequence of elements in
a positively graded polynomial ring $P$, generating an ideal $I$.
Suppose that the leading forms $\bar f_{1},\ldots,\bar f_{r}$ form a
regular sequence.  Then the $f_{i}$ form a regular sequence, and the
ideal of leading forms of $I$ is generated by the $\bar f_{i}$.\qed
\end{lemma}

We conclude that the associated graded ring of the $v$--filtration of
$X(\Delta)$ is a
two-dimensional
complete intersection, defined by the $t-2$ leading forms as above.
In addition, $X(\Delta)$ is a complete intersection at the origin.
\end{proof}

\begin{corollary}\label{isolated} Let $X(\Delta)$ be a splice diagram singularity as
above.  Then for any two variables $z,z'$, the locus
$$X(\Delta)\cap\{z=z'=0\}$$
contains the origin as an isolated point.
\end{corollary}
\begin{proof}  Choose $v$ the node adjacent to the leaf
corresponding to $z$.  In Spec of the associated graded, the locus
$\{z=z'=0\}$ is exactly the origin, so the associated graded modulo
$z$ and $z'$ is a 0-dimensional complete intersection.  Thus, the
local ring of $X$ at the origin modulo these two variables is also
0-dimensional.
\end{proof}

One can also see that any of these associated graded rings is reduced;
it suffices to show generic reducedness.  Choose a point all of whose
coordinates are non-zero.  Then in every sector, one can solve the
splice equations around that point and set each variable $z$ equal to
a constant times a power of a new variable $t$.  This transforms that
sector's admissible monomial into a power of one variable $t$.  This
puts the equations corresponding to $v$ into the form of a Brieskorn
complete intersection.  In particular, Spec of the associated graded
is smooth at points for which all coordinates are non-0; the only
singular points occur along the curves obtained by setting a variable
equal to 0.

To show the singularity of $X(\Delta)$ is isolated, and to get a
handle on resolution diagrams, we do a weighted blow-up.

\section{Weighted blow-up and the proper transform}\label{sec:blow-up}

Let $z_{1},\dots,z_t$ be coordinates on an affine space $\C
^{t}$, where $z_{i}$ has positive integer weight $m_{i}$.
Blowing-up
the corresponding weight filtration gives the \emph{weighted blow-up}
$Z\rightarrow \C^t$, an isomorphism off the inverse image of
the origin.  $Z$ is covered by $t$ affine varieties $U_{i}$,
each of
which is a quotient of an affine space $V_{i}$ by a cyclic group of order
$m_{i}$.  $V_{1}$ has coordinates $A_{1},\dots,A_{t}$,
 related to the $z_{i}$ via
$$z_{1}=A_{1}^{m_{1}},~ z_{2}=A_{1}^{m_{2}}A_{2},~ \dots,~
z_{t}=A_{1}^{m_{t}}A_{t};$$ $U_{1}$ equals $V_{1}$
modulo the action of the
cyclic
group generated by $$S=[-1/m_{1},m_{2}/m_{1},\dots,m_{t}/m_{1}]\,,$$
where we are using the notation
$$[q_1,\dots,q_t]:=(\exp(2\pi i q_1),\dots,\exp(2\pi iq_t))\,.$$
Consider now a splice diagram singularity $X=X(\Delta)$, where as
usual $t$ is the number of leaves.  Let $v^{*}$ be an \emph{end-node}
of $\Delta$, ie, a node all but one of whose edges is adjacent to a
leaf.  Putting this node on the left side of $\Delta$, we write
$$\xymatrix@R=6pt@C=36pt@M=0pt@W=0pt@H=0pt{
\Circ&&&&\\
\Vdots&\undertag\Circ{v*}{2pt}\lineto[ul]_(.35){n_1}
\lineto[dl]^(.35){n_\sigma}\lineto[rr]^(.2)r
&&\Circ\dashto[ur]\dashto[dr]&\Vdots\\
\Circ&&&&}
$$
(As the one-node case is easy and is revisited below, we assume there
are at least two nodes).  Blow-up the corresponding weight filtration,
where $z_{w}$ has weight $\ell_{v^{*}w}$.  Since the associated graded
of $X$ with respect to the filtration is a complete intersection, the
proper transform of $X$ in $Z$ or the $V_{i}$ is defined by the proper
transform of the splice equations.  Further, since the origin is the
only point of $X$ at which two coordinates vanish, one sees the full
transform of $X$ by considering the proper transforms $X_{i}$ on
$V_{i}$ for $i=1,2$.

Denote by
$\sigma$ the number of leaves adjacent to $v^{*}$, $\tau=t-\sigma$ the
number of remaining leaves and number these leaves
$w_1,\dots,w_{\sigma+\tau}$.  Call the corresponding variables
$x_{1},\dots,x_{\sigma},y_{1},\dots,y_{\tau}$.  Write
$N=n_{1}n_{2}\dots n_{\sigma}, N_{i}=N/n_{i}$.  Then the weight of
$x_{i}$ is $rN_{i}$, while that of $y_{k}$ is $N\ell'_k$, where
$\ell'_k:=\ell'_{v^*w_{\sigma+k}}$.  We write the splice diagram
equations in two groups.

First, one has
\begin{subequations}\label{eq:general splice}
\begin{equation}\label{eq:13}
x_{i}^{n_{i}}+a_{i}x_{\sigma}^{n_{\sigma}}=~
b_{i}\underline{y}^{\underline{\alpha}}+ H_{i}(\underline x,\underline
y),\quad 1\leq i\leq \sigma -1\,.
\end{equation}
Here, $\underline{y}^{\underline{\alpha}}:=
\prod_{k=1}^{\tau}y_{k}^{\alpha_{k}}$ is an admissible monomial for
$v^{*}$, so $\sum \alpha_{k}\ell'_{{k}}=r$.  The term $H_i(\underline
x,\underline y)$ (shorthand for
$H_i(x_1,\dots,x_\sigma,y_1,\dots,y_\tau)$) contains monomials of
higher $v^{*}$--weight.

Second, for each node $v\neq v^{*}$, order the
admissible monomials $M_{vj},j=1,\dots,\delta_{v}$ so that $M_{v1}$
corresponds to the edge pointing in the direction of $v^{*}$ (thus
only $M_{v1}$ can involve any $x_{i}$ variables). The corresponding
equations are of the form
  \begin{equation}
  \label{eq:14}
  \sum_{j=1}^{\delta_{v}} a_{vij}M_{vj}=H_{vi}(\underline x,\underline
y)\,,\quad
i=1,\dots,\delta_{v}-2\,.
\end{equation}
\end{subequations}
The transforms of the first group of equations (\ref{eq:13})
on $V_{1}$ have the form
\begin{subequations}\label{eq:transformed general splice}
\begin{equation}
  \label{eq:8}
  \begin{split}
1+a_{1}A_{\sigma}^{n_{\sigma}}&=b_{1}\underline{\z}^{\underline{\alpha}}+
A_{1}J_{1}(\underline A,\underline \z)\\
A_{i}^{n_{i}}+a_{i}A_{\sigma}^{n_{\sigma}}&=
b_{i}\underline{\z}^{\underline{\alpha}}+
A_{1}J_{i}(\underline A,\underline \z)\,,\quad 2\leq i\leq \sigma -1\,.
\end{split}
\end{equation}
For equations (\ref{eq:14}) corresponding to another node $v$, it follows from
Lemma \ref{le:weights}
that the transforms of such
equations take the form
\begin{equation}
  \label{eq:9}
  \sum_{j=2}^{\delta_{v}}
a_{vij}M_{vj}(\underline \z)=A_{1}H'_{vi}(\underline A,\underline \z),\quad
i=1,\dots,\delta_{v}-2\,.
\end{equation}
\end{subequations}
Here, $M_{vj}$ (for $j>1$) is the same monomial as before, but
evaluated on the $\z_{k}$ instead of the $y_{k}$.  These equations are
obtained from the ones above by substituting for the $x_{i}$ and
$y_{k}$, and then dividing through by the highest power of $A_{1}$
that occurs, which is $A_1^{\ell_{v^{*}v}}$ by Lemma \ref{le:weights}.

One needs to get a handle on the singular locus of $X_{1}$.  First, the
exceptional divisor is the subscheme of $X_{1}$ defined by $A_{1}=0$;
it is
given by equations
\begin{align*}
1+a_{1}A_{\sigma}^{n_{\sigma}}&=b_{1}\underline{\z}^{\underline{\alpha}}\,,\\
A_{i}^{n_{i}}+a_{i}A_{\sigma}^{n_{\sigma}}&=
b_{i}\underline{\z}^{\underline{\alpha}}\,,&2\leq i\leq \sigma -1\,,\\
\sum_{j=2}^{\delta_{v}}
a_{vij}M_{vj}(\underline \z)&=0\,,
&v\ne v^* \text{ a node,}\quad i=1,\dots,\delta_{v}-2\,.
\end{align*}

\begin{lemma} \label{le:4.1}
The curve defined above has singularities only at the $n_2\dots
n_\sigma $ points for which all $\z_{k}=0$ (hence no $A_{i}=0$ for
$i>1$).  Every connected component of the curve contains such a
point.  At a point for which $A_{2}=0$, $A_{2}$ is a local
analytic coordinate; there are $n_3\dots n_\sigma r$ such
points.
\end{lemma}
\begin{proof}
 The equations arising from
  $v$ different from $v^{*}$ are splice diagram curve equations for
  the splice diagram obtained by removing from $\Delta$ the $\sigma$
  edges and leaves adjacent to $v^{*}$, which is now viewed as a
  root.   By Theorem \ref{th:curves},
  these define a reduced complete intersection
  curve $C$ in the variables $\z_{k}$, with one singularity at the origin.
  Adding in the first $\sigma -1$ equations (and variables
  $A_{2},\dots,A_{\sigma}$) defines a branched cover of this curve,
  which by the Jacobian criterion is unramified except when one of the
  $A_{i}$ is 0.  There are thus $N_{1}=n_{2}\dots n_{\sigma}$
  singular points lying above the origin.  The image of any connected component
  under the finite cover
  must contain the origin, whence each component contains at least
  one of the $N_{1}$ points.

  Next consider a point where $A_{2}=0$.  The genericity condition on
  the coefficients implies all other $A_{i}\neq 0$.  Some $\z_{j}\neq
  0$, hence all $\z_{k}\neq 0$, by Theorem \ref{th:curves}.  In a
  neighborhood of such a point, we can write
  $\z_{i}=c_{i}t^{\ell_{i}}$, as in the proof of Theorem
  \ref{th:curves}, and  so replace the term
  $\underline{\z}^{\underline{\alpha}}$ by a constant times $t^r$.  So
  the curve is now defined by $\sigma -1$ equations in the variables
  $A_{2},\dots,A_{\sigma},t$.  Again the Jacobian criterion implies
  that $A_{2}$ is a coordinate at any point where all the other
  coordinates are non-zero (one needs again the precise genericity of
  the coefficients). To count the points we note that the value of
  $\underline \z^{\underline \alpha}$ is determined, so part 5 of
  Theorem \ref{th:curves} gives us exactly $r$ points on the curve in
  the variables $\z_k$, and above each of these there are $n_3\dots
  n_\sigma$ points when one adds in the variables $A_2=0,
  A_3,\dots,A_\sigma$.
\end{proof}

Thus, the singular locus of $X_{1}$ intersects the exceptional divisor
only at the $N_{1}$ singular points described above.  We show these
are isolated singular points of $X_{1}$.  Repeating for $X_{2}$ will
imply that $X$ has an isolated singularity.  In addition, we will have the
necessary set-up to study singular points on $X_{1}$ and on its quotient by
finite groups (as needed in Section \ref{sec:induct}).

Choose a point where $A_{1}$ and all the $\z_{k}$ are 0; thus, all
$A_{i}$ are non-zero for $i>1$.  Now, the matrix of $A_{i}$--partial
derivatives ($2\leq i \leq \sigma $) of the $\sigma -1$ equations
(\ref{eq:8}) is easily seen to be invertible at such a point.  So, by
the implicit function
theorem one can, in a neighborhood of such a point, solve uniquely
these equations, and write each $A_{i}$ as a convergent power series
in $A_{1}, \z_{1},\dots,\z_{\tau}$, with non-zero constant term.  Plug
these convergent power series into the second group of equations
(\ref{eq:9}).  We show one now has a set of splice diagram equations
for a smaller diagram $\tilde \Delta$, which by induction represents
an isolated singularity.  This will complete the proof of Theorem
\ref{th:splice is CI}.

Let $\tilde\Delta$ be
the splice diagram whose underlying graph is $\Delta$ less the
$\sigma$ edges and leaves adjacent to $v^{*}$; thus $v^{*}$ is now
replaced by a leaf $w^{*}$.  Edge-weights not pointing towards $v^{*}$
are defined to be the same as before.  For a node $v$ and the edge at
$v$ pointing towards $v^{*}$, define a new edge-weight by
\begin{equation}
  \label{eq:11}
\tilde
 d_{v1}=rd_{v1}-N(d_{v}/d_{v1})(\ell'_{vv^{*}})^{2}
\end{equation}
 This is an integer, and its positivity is easily seen by
 multiplying all the edge-determinant inequalities between $v^{*}$ and
 $v$.  One readily checks by induction over distance of an edge from
 $v^*$ that each edge determinant of $\tilde\Delta$
 is $r$ times the corresponding edge-determinant of $\Delta$, hence
 positive.

 \begin{lemma}\label{le:tilde}
 Assign the variables $A_{1}$ to $w^{*}$ and $\z_{k}$ to
 corresponding other leaves of $\tilde\Delta$.  Then the proper transforms
 of the second group of equations for $\Delta$, with substitutions
 for $A_{2},\dots,A_{\sigma}$, are (in a neighborhood of the point
 in question) splice
 diagram equations for $\tilde\Delta$ (which in particular must satisfy
 the semigroup condition).
 \end{lemma}
 \begin{proof}  Choose a node $v$ of $\tilde\Delta$.
   Since the edge-weights of $\tilde\Delta$ are the same as the
   corresponding ones of $\Delta$ except on edges pointing to $w^{*}$,
   one sees that the old admissible monomials $M_{vi}(\underline y)$
   for $i>1$ and $\Delta$ yield admissible monomials
   $M_{vi}(\underline \z)$ for $\tilde\Delta$.  One must check the
   replacement for $M_{v1}(\underline x,\underline y)$ in the new
   equations.

 More generally, consider any monomial
 $$\prod_{i} x_{i}^{\beta_{i}}\prod_{k}y_{k}^{\gamma_{k}}$$
 appearing in an equation associated to the node $v$.
 Going up to $X_{1}$ means substituting for
 $x_{i}$ and $y_{k}$ in terms of the $A_{i}$ and $\z_{k}$; taking
 proper transform means subtracting $\ell_{v^{*}v}$ from the exponent of
 $A_{1}$.  This gives a monomial of the form
 $$A_{1}^{Q}A_{2}^{\beta_{2}}\dots A_{\sigma}^{\beta_{\sigma}}\prod
 \z_{k}^{\gamma_{k}}\,.$$
 Here, one has $$Q~=~r\sum \beta_{i}N_{i}~+~N\sum
 \gamma_{k}\ell'_{{k}}~-~\ell_{vv^{*}}\,.$$
 In a neighborhood of the singular point,
 $A_{i}$ for $i>1$ is a power series in the variables
 $A_{1},\z_{k}$ with non-zero constant term; so we are really
 considering (up to a fixed factor plus higher-order terms) the
 transformed monomial
 $$A_{1}^{Q}\prod \z_{k}^{\gamma_{k}}\,.$$
 If the original monomial is the particular admissible $M_{v1}$ for
 $\Delta$, one has
 \begin{equation}
   \label{eq:2}
   \sum_{i=1}^{\sigma}\beta_{i}\ell_{vw_i}+\sum_{k=1}^\tau
 \gamma_{k}\ell_{vw_{\sigma+k}} ~=~d_{v}\,,
 \end{equation}
 the second sum being over leaves on the $v^{*}$ side of $v$.   Let us
 put a $~\tilde{}~$ over a
 linking number or degree computed in $\tilde\Delta$.  To prove the
 transformed monomial is admissible for $\tilde\Delta$, one must prove
 \begin{equation}
   \label{eq:3}
Q\tilde{\ell}_{vv^{*}}+\sum_{k=1}^\tau
\gamma_{k}\tilde{\ell}_{vw_{\sigma+k}}~
 =~\tilde{d}_{v}\,.
\end{equation}
We claim that, in fact,
\begin{equation}\label{eq:4}
Q\tilde{\ell}_{vv^{*}}+\sum_{k=1}^{\tau}
\gamma_{k}\tilde{\ell}_{vw_{\sigma+k}}~
 -\tilde{d}_{v}=r\left(
\sum_{i=1}^{\sigma}\beta_i\ell_{vw_i}+
\sum_{k=1}^{\tau}\gamma_k\ell_{vw_{\sigma+k}}-
d_{v}\right)\,,
\end{equation}
so (\ref{eq:3}) is equivalent to (\ref{eq:2}).  We postpone the proof
of (\ref{eq:4}), which holds even if some of the leaves $w_{\sigma+k}$ are
beyond $v$ from $v^*$ (we need this later).

 We conclude that our equations have appropriate
 admissible monomials; in particular, $\tilde\Delta$ satisfies
 the semigroup condition.  Further, the coefficients at each node
 satisfy the appropriate genericity condition, since they are the
 same as before except for a multiple of the fixed factor just
 mentioned.

 Next, one must show a monomial of $v$--weight $>d_{v}$ gives rise to a
 transformed monomial of weight $>\tilde{d_{v}}$ with respect to the
 new $\tilde{v}$--valuation.  Thus, for non-negative integers
 $\beta_{i}$, $\gamma_{k}$, the inequality
 $$\sum_{i=1}^{\sigma}\beta_{i}\ell_{vw_i} +\sum_{k=1}^{\tau}
 \gamma_{k}\ell_{vw_{\sigma+k}}~>~d_{v}$$
 should imply the inequality
 $$Q\tilde{\ell}_{vv^{*}} +\sum_{k=1}^{\tau}
\gamma_{k}\tilde{\ell}_{vw_{\sigma+k}} ~>~ \tilde{d}_{v}\,.$$
The equivalence of these inequalities also follows from equation
 (\ref{eq:4}).

 To complete the proof we must thus prove equation (\ref{eq:4}).
 First note that multiplying equation (\ref{eq:11}) by
 $d_v/d_{v1}$ gives
\begin{equation*}
\tilde d_v=rd_v-N(\frac{d_v}{d_{v1}}\ell'_{vv^*})^2=
rd_v-\tilde\ell_{vv^*}\ell_{vv^*}\,.
\end{equation*}
Hence
\begin{align*}
  &Q\tilde{\ell}_{vv^{*}}+\sum \gamma_{k}\tilde{\ell}_{vw_{\sigma+k}}~
  -\tilde{d}_{v}\\
  =~&\left(r\sum \beta_{i}N_{i}\ +N\sum \gamma_{k}\ell'_{{k}}~
    -\ell_{vv^{*}}\right)\tilde\ell_{vv^*}+ \sum
  \gamma_{k}\tilde{\ell}_{vw_{\sigma+k}}~
  -\left(rd_v-\tilde\ell_{vv^*}\ell_{vv^*}\right) \\
  =~&r\sum\beta_iN_i\tilde\ell_{vv^*}
  +\sum\gamma_k\left(N\ell'_{{k}}\tilde\ell_{vv^*}+
    \tilde\ell_{vw_{\sigma+k}}\right) - rd_{v}\\
  =~&r\sum\beta_i\ell_{vw_i}+r\sum\gamma_k\ell_{vw_{\sigma+k}} - rd_{v}\,,
\end{align*}
where the equality $N\ell'_{{k}}\tilde\ell_{vv^*}+
\tilde\ell_{vw_{\sigma+k}}=r\ell_{vw_{\sigma+k}}$ is seen as follows: Denote by $v_k$
the vertex where the paths from $v$ to $v^*$ and $w_{\sigma+k}$ diverge (so
$v_k=v$ iff $w_{\sigma+k}$ is beyond $v$ from $v^*$). Then
\begin{gather*}
  N\ell'_{{k}}\tilde\ell_{vv^*}+
  \tilde\ell_{vw_{\sigma+k}}=\ell'_{{k}}\ell_{vv^*}+
  \frac{\ell_{vw_{\sigma+k}}}{d_{v_k}}\tilde d_{v_k}=
  \ell'_{{k}}\ell_{vv^*}+
  \frac{\ell_{vw_{\sigma+k}}}{d_{v_k}}\left(rd_{v_k}-\tilde\ell_{v_kv^*}
    \ell_{v_kv^*}\right)\\
  =r\ell_{vw_{\sigma+k}}+\ell'_{{k}}\ell_{vv^*}-
  \ell_{vw_{\sigma+k}}\tilde\ell_{v_kv^*}\ell_{v_kv^*}/d_{v_k}=
  r\ell_{vw_{\sigma+k}}\,,
\end{gather*}
since it is easy to check that
$\ell'_{{k}}\ell_{vv^*}=\ell'_{v^*w_{\sigma+k}}\ell_{vv^*}$ and
$\ell_{vw_{\sigma+k}}\tilde\ell_{v_kv^*}\ell_{v_kv^*}/d_{v_k}$
represent the same product of weights.  This completes the proof of
Lemma \ref{le:tilde}, and hence also of Theorem \ref{th:splice is CI}.
\end{proof}

\section{The discriminant group and its natural representation}\label{sec:discriminant}

We consider the dual resolution graph $\Gamma$ of a good resolution of
a normal surface  singularity with rational homology sphere link
(definitions relating to dual resolution graph, minimal good
resolution, etc., are recalled in Section \ref{sec:splicing}).  A
vertex $v$ in $\Gamma$ corresponds to an exceptional curve $E_{v}$,
and an edge corresponds to an intersection of two exceptional curves.
A vertex $v$ is called a \emph{leaf} (or \emph{end}) if its valency is
1, a \emph{node} if its valency is $\geq 3$. Each vertex $v$
is weighted by self-intersection number of its associated curve $E_v$.

Let $$\E :=
\bigoplus_{v\in vert(\Gamma)}~ \Z \cdot E_{v}$$
be the lattice
generated by the classes of these curves (so $\E$ can be identified
with the homology $H_{2}(\bar X;\Z)$ of the resolution).  Via the
negative-definite intersection pairing $A(\Gamma)$, one has natural
inclusions $$\E \subset \E ^{\star} =\Hom(\E , \Z) \subset \E \otimes
\Q\,.$$
The discriminant group is the finite abelian group
$$D(\Gamma):= \E ^{\star}/\E\,,$$
whose order is
$\disc(\Gamma):=\det(-A(\Gamma))$.  There are induced symmetric pairings of
$\E \otimes \Q$ into $\Q$ and $D(\Gamma)$ into $\Q /\Z$.

To calculate the discriminant group, let $\{e_{v}\}\subset \E
^{\star}$ be the dual basis of the ${E_{v}}$, ie,
$$e_{v}(E_{v'})=\delta_{vv'}\,.$$
We claim that the images of those
$e_{w}$ for which $w$ is a leaf of the graph generate $D(\Gamma)$.  In
fact, more is true:

\begin{proposition} Consider a collection ${e_{w}}$, where $w$ runs
  through all but one leaf of the graph $\Gamma$.  Then $D(\Gamma)$ is
  generated by the images of these $e_{w}$.
\end{proposition}
\begin{proof}
  Let $E_{v}$ be any exceptional curve, with neighbors
  $E_{1},\dots,E_{r}$ ($r\geq 1$).  Then in $\E ^{\star}$, one
  verifies by dotting with any curve that
  $$
  E_{v}=(E_{v}\cdot E_{v})e_{v}+\sum_{i=1}^{r}e_{i}.
  $$
  Thus, in $D(\Gamma)$ any $e_{i}$ can be completely expressed in
  terms of the $e_{v}$ corresponding to one neighbor and all the other
  curves on the far side of that neighbor.  Choose any one end curve
  $E_{w}$; then every curve $E_{v}$, except for the remaining end
  curves, has a neighbor away from $E_{w}$, and so the corresponding
  $e_{v}$ may be expressed in terms of outer curves.  Eventually, all
  are expressed in terms of the remaining end curves.
\end{proof}

\begin{proposition}\label{prop:free off hyperplanes}
Let $e_{1},\dots,e_{t}$ be the elements of the
  dual basis of\/ $\E^{\star}$ corresponding to the $t$ leaves of
  $\Gamma$.  Then the homomorphism $\E^{\star}\rightarrow \Q^{t}$
  defined by
  $$e\mapsto (e\cdot e_{1},\dots,e\cdot e_{t})
  $$
  induces an injection
  $$D(\Gamma)=\E ^{\star}/\E ~\hookrightarrow~(\Q /\Z )^{t}.
  $$
  In fact, each non-trivial element of $D(\Gamma)$ gives an element
  of $(\Q /\Z )^{t}$ with at least two non-zero entries.
\end{proposition}
\begin{proof}
  It suffices to show that if $e\cdot e_{i} \in \Z$ for $1\leq i \leq
  t-1$, then $e\in \E$. But then the set of $e'\in \E^{\star}$ for
  which $e\cdot e' \in \Z$ is a subgroup containing $\E$ and these
  $e_{i}$, so by the last proposition must be all of $\E ^{\star}$.
  Write $e=\sum r_{v}E_{v}$ as a rational combination of exceptional
  curves; then for every exceptional curve $E_{v}$ one has $e\cdot
  e_{v}=r_{v}~\in \Z$, as desired.
\end{proof}

It will be convenient to exponentiate, and to consider
$$(\Q /\Z)^{t} \hookrightarrow (\C ^{\star})^{t}$$
via $$(\dots,
r,\dots) \mapsto (\dots, \exp(2\pi ir),\dots)=:
[\dots,r,\dots]\,.$$
%Thus for a leaf $w$ the discriminant action on $\C^t$ is via
% $$e_{w}\mapsto [e_{w}\cdot e_{w_1},~\dots~,e_{w}\cdot e_{w_t}].$$
Keeping in mind the last proposition, we
summarize in the
\begin{proposition} \label{prop:5.3}
  If the leaves of\/ $\Gamma$ are numbered $w_1,\dots,w_t$, then the
  discriminant group $D(\Gamma)$ is naturally represented by a
  diagonal action on $\C ^{t}$, where the entries are $t$--tuples of
  $\disc(\Gamma)$-th roots of unity.  Each leaf\/ $w_j$ corresponds to
  an element
  $$[e_{w_j}\cdot e_{w_1},\dots,e_{w_j}\cdot e_{w_t}]:=\left(\exp(2\pi
    i\,e_{w_j}\cdot e_{w_1}),\dots,\exp(2\pi i\,e_{w_j}\cdot
    e_{w_t})\right)\,,$$
  and any $t-1$ of these generate $D(\Gamma)$.
  The representation contains no pseudo-reflections, ie,
  non-identity elements fixing a hyperplane.
\end{proposition}

\section{Resolution graphs and the congruence condition}
\label{sec:congruence condition}

Non-minimal resolutions are needed later in our inductive arguments, so
we do not want to insist that $\Gamma$ corresponds to the minimal good
resolution. We therefore make the following purely technical definition.

\begin{definition}[\Minimal{}ity]\label{minimality}
  The resolution tree $\Gamma$ is \emph{\minimal} if any string in
  $\Gamma$ either contains no $(-1)$--weighted vertex, or consists
  of a unique $(-1)$--weighted vertex (a \emph{string} is a
connected subgraph that includes no node of $\Gamma$).
\end{definition}

Associated to a string\comment{06/25 Added. Being careful not to have
  to assume \minimal{}ity in first part of Lemma \ref{le:cf}}
$$\xymatrix@R=6pt@C=24pt@M=0pt@W=0pt@H=0pt{
  E=&\overtag{\Circ}{-b_1}{8pt}\lineto[r]&
  \overtag{\Circ}{-b_2}{8pt}\dashto[r]&
  \dashto[r]&\overtag{\Circ}{-b_k}{8pt}\\&~}$$
in a resolution diagram is a continued fraction
    $$n/p=b_{1}-1/b_{2}-1/\dots -1/b_{k}\,.$$
The continued fraction $1/0$ is associated with the empty string.
We will need the following standard facts about this relationship,
whose proofs are left to the reader.
\begin{lemma}\label{le:cf}
    Reversing a string with continued fraction $n/p$ gives one with
  continued fraction $n/p'$ with $pp'\equiv 1$ $($mod $n)$. Moreover,
  the following hold:
  \begin{align*}
n&=\disc\left(\hspace{8pt} \xymatrix@R=6pt@C=24pt@M=0pt@W=0pt@H=0pt{
    \overtag{\Circ}{-b_1}{8pt}\lineto[r]&
  \overtag{\Circ}{-b_2}{8pt}\dashto[r]&
  \dashto[r]&\overtag{\Circ}{-b_k}{8pt}}\hspace{8pt} \right)\\
p&=\disc\left(\hspace{8pt} \xymatrix@R=6pt@C=24pt@M=0pt@W=0pt@H=0pt{
    \overtag{\Circ}{-b_2}{8pt}\dashto[r]&
  \dashto[r]&\overtag{\Circ}{-b_k}{8pt}}\hspace{8pt} \right)\\
p'&=\disc\left(\hspace{8pt} \xymatrix@R=6pt@C=24pt@M=0pt@W=0pt@H=0pt{
    \overtag{\Circ}{-b_1}{8pt}\dashto[r]&
  \dashto[r]&\overtag{\Circ}{-b_{k-1}}{8pt}}\hspace{8pt} \right)\,,
\end{align*}
and the continued fraction in the last case is $p'/n'$ with
$n'=(pp'-1)/n$.

There is a unique directed \minimal{} string for each
$n/p\in[1,\infty]$, and in this case the reversed string has continued
fraction $n/p'$ with $p'$ the unique $p'\le n$ with $pp'\equiv 1$ $($mod
$n)$.\qed
\end{lemma}

Associate to a (not necessarily minimal)
resolution graph $\Gamma$ a splice diagram $\Delta$,
as in \cite{neumann-wahl02} (see also section \ref{sec:splicing}): First,
suppress all vertices of valency two in $\Gamma$, yielding a tree of
the same general shape, but now with only leaves and nodes. Second, to
every node $v$ and adjacent edge $e$ of $\Delta$ (or $\Gamma$),
associate an edge-weight $d_{ve}$ as follows: removing the node from
$\Gamma$, take the positive determinant $d_{ve}:=\disc(\Gamma_{ve})$
of the remaining connected graph $\Gamma_{ve}$ in the direction of the
edge.  The splice diagram has positive edge-determinants and satisfies
the ideal condition (Definition \ref{def:ideal}; this is proved in
section \ref{sec:splicing}).  In the unimodular case, with
$\disc(\Gamma)=1$, the weights around a node are relatively prime; but
this is no longer true in general.

The discriminant group $D(\Gamma$) acts diagonally on $\C ^{t}$, as in
Proposition %\ref{prop:free off hyperplanes}.
\ref{prop:5.3}. Viewing the $z_{w}$'s as
linear functions on $\C ^{t}$, $D(\Gamma)$ acts naturally on the
polynomial ring $P=\C [\dots,z_{w},\dots]$; $e$ acts on monomials as
$$\Pi z_{w}^{\alpha_{w}} \mapsto \left[-\sum (e\cdot e_{w})\alpha_{w}\right]
\Pi z_{w}^{\alpha_{w}}\,.$$  In other words, the group transforms this
monomial by multiplying by the character
$$e \mapsto \exp\left(-2\pi i\sum (e\cdot e_{w})\alpha_{w}\right)\,.$$
If $\Delta$ satisfies the semigroup condition, one has the notion of
admissible monomials (Definition \ref{admissible}).
\begin{definition}[Congruence Condition]\label{def:congruence}
  Let $\Gamma$ be a resolution diagram, yielding a splice diagram
  $\Delta$ satisfying the semigroup condition.  We say $\Gamma$
  satisfies the \emph{congruence condition} if for each node $v$, one
  can choose for every adjacent edge $e$ an admissible monomial
  $M_{ve}$ so that $D(\Gamma)$ transforms each of these monomials
  according to the same character.
\end{definition}

We can write down this condition explicitly in terms of $\Gamma$
 and $\Delta$.
\begin{lemma}\label{le:disc is inverse}
The matrix $(e_{v}\cdot e_{v'})$ $(v,v'$ vertices of
$\Gamma)$ is the inverse
of the matrix $A(\Gamma)=(E_{v}\cdot E_{v'})$.
\end{lemma}
\begin{proof}  By elementary linear algebra, the matrix of the dual basis
  in an inner product space (such as $\E \otimes \Q$) is the inverse
  of the matrix of the original basis.
\end{proof}

\begin{proposition}\label{prop:ell}
Let $w,w'$ be distinct leaves of\/ $\Gamma$,
  corresponding to distinct leaves of $\Delta$, and let $\ell _{ww'}$
  denote their linking number.  Then
  $$
  e_{w}\cdot e_{w'}=-\ell _{ww'}/\disc(\Gamma).
  $$
\end{proposition}
\begin{proof}
This is immediate by the preceding lemma and Theorem \ref{th:props}.
\end{proof}

This proposition implies that for distinct leaves $w,w'$ the number
$e_w\cdot e_{w'}$ depends (except for the denominator $\disc(\Gamma)$)
only on the splice diagram $\Delta$. For a leaf $w$, the number $(e_{w}\cdot
e_{w})\disc(\Gamma)$ is not determined solely by $\Delta$.

\begin{proposition}\label{prop:endwt}\comment{06/25: Changed, there
    are no minimality assumptions now; proofs is slightly simpler}
  Suppose we have a string from a node $v$ to an adjacent leaf $w_1$
  in $\Gamma$ with associated continued fraction $d_1/p$,
  so $d_1$ is the weight at $v$ towards $w_1$. Let $p'$ be the
determinant  %discriminant
of the same string with the last vertex $w_1$ deleted,
  so
  $pp'\equiv 1$ $($mod $d_{1})$ (see Lemma \ref{le:cf}).  Then (with
  $d_v$ the product of weights at $v$)
$$e_{w_1}\cdot e_{w_1}=-d_{v}/\left(d_{1}^{2}\disc(\Gamma)\right) ~-~
p'/d_{1}\,.$$
(Compare this with $e_{w_1}\cdot
e_{w_2}=-d_v/\left(d_1d_2\disc(\Gamma)\right)$ for two leaves adjacent
to $v$.)
\end{proposition}
\begin{proof}  Lemma \ref{le:disc is inverse}, Theorem \ref{th:props},
  and Lemma \ref{le:end weight}.
\end{proof}
\begin{corollary}\label{cor:action} The class of $e_{w'}$ ($w'$ a leaf) transforms the
  monomial $\Pi z_{w}^{\alpha_{w}}$ by multiplication by the root of
  unity $$\Bigl[\sum_{w \neq w'}\alpha _{w}\ell_{ww'}/\disc(\Gamma)~-~
  \alpha_{w'}e_{w'}\cdot e_{w'}\Bigr].\eqno{\qed}$$
\end{corollary}

These formulas allow a direct way to check the congruence condition.

\begin{proposition}\label{prop:congruence}
 Let $\Gamma$ be a graph whose splice diagram
  $\Delta$ satisfies the semigroup condition.  Then the congruence
  condition is equivalent to the following: for every node $v$ and
  adjacent edge $e$, there is an admissible monomial $M_{ve}=\Pi
  z_{w}^{\alpha_{w}}$ ($w$ running through the leaves in
  $\Delta_{ve}$) so that for every leaf $w'$ of $\Delta_{ve}$,
$$\Bigl[\sum_{w \neq w'}\alpha _{w}\ell_{ww'}/\disc(\Gamma)~-~
\alpha_{w'}e_{w'}\cdot e_{w'}\Bigr]=\left[\,\ell_{vw'}/\disc(\Gamma)\,\right]\,.$$
\end{proposition}
\begin{proof} We first
claim that if $\bar{e}\neq e$ is another edge of $v$,
then $e_{w'}$ transforms any admissible polynomial $M_{v\bar{e}}$ by
the root of unity
$$[\,\ell_{vw'}/\disc(\Gamma)\,]\,.$$
To see this, one checks (via the definition of linking numbers) that for
$\bar{w}\in \Delta_{v\bar{e}}$, one has
$$\ell_{w'\bar{w}}=\ell_{w'v} \ell'_{v\bar{w}}/d_{v\bar{e}}\,.$$
In particular, if $$d_{v\bar{e}}=\sum
\beta_{\bar{w}}\ell'_{v\bar{w}}\,,$$ then $$\sum
\beta_{\bar{w}}\ell'_{w'\bar{w}}=\ell_{vw'}\,.$$  The claim follows.

In particular, at each node $v$, checking the congruence condition on
the $M_{ve}$'s imposed by one $e_{w}$ involves only the stated
equality, involving the edge in the direction of $w$.
\end{proof}

Note that there is nothing to check for edges leading to leaves.  In
case $\Gamma$ is star-shaped (ie, $\Delta$ has only one node), there
are no semigroup conditions, hence no congruence conditions.  We will
later explain these conditions completely in the two-node case.  But
we give one example.

\begin{example}{\rm\cite{neumann-wahl02}}\qua Consider the resolution diagram
$$  \xymatrix@R=6pt@C=24pt@M=0pt@W=0pt@H=0pt{
\\
&\overtag{\Circ}{-3}{8pt}&&&\overtag{\Circ}{-3}{8pt}\\
{\Gamma=}&&\overtag{\Circ}{-7}{8pt}\lineto[ul]\lineto[dl]\lineto[r]
&\overtag{\Circ}{-1}{8pt}\lineto[ur]\lineto[dr]\lineto[l]&\\
&\overtag{\Circ}{-3}{8pt}&&&\overtag{\Circ}{-3}{8pt}\\&~}\,.$$
The corresponding splice diagram is
$$
\xymatrix@R=6pt@C=24pt@M=0pt@W=0pt@H=0pt{
\\
&\Circ&&&&\Circ\\
{\Delta=}&&\Circ\lineto[ul]_(.25){3}\lineto[dl]^(.25){3}\lineto[rr]^(.25){3}^(.75){57}&
&\Circ\lineto[ur]^(.25){3}\lineto[dr]_(.25){3}&\\
&{\Circ}&&&&{\Circ}\\&~\\&~}
$$
Labeling the variables $x,y,u,v$ clockwise from the bottom left
leaf,\comment{06/24: This should be an example for the proposition, so
  I have changed to do it that way} one checks that an admissible monomial
$u^\alpha v^\beta$ at the left node satisfies the condition of
Proposition \ref{prop:congruence} if and only if $\alpha$ and $\beta$
are both $\equiv 2$ (mod $3$). This is incompatible with the
admissibility condition $\alpha+\beta=1$, so $\Gamma$ does not satisfy
the congruence condition.\comment{We could leave in: ``Alternatively,
  admissible monomials for the left node are $x^3,y^3$, and $u$ or
  $v$.  One verifies that $\disc(\Gamma)=90$, and that the element of
  the discriminant group corresponding to the bottom right leaf gives
  the action on the variables given by $T=[1/10,1/10,19/30,29/30]$. We
  see that $x^3$ and $y^3$ do not transform as do either $u$ or $v$.''}
\end{example}

\section{Splice diagram equations with discriminant group action}

Let $\Gamma$ be a \minimal{} resolution tree
(Definition \ref{minimality}).
Assume $\Gamma$ satisfies the semigroup and congruence conditions
(\ref{def:sg}, \ref{def:congruence}).
Let $\Delta$ be the corresponding splice diagram, and $z_{w}$ a
variable associated to each leaf $w$.  The discriminant group
$D(\Gamma)$ acts on the monomials in the variables $z_{w}$.

For each node $v$, choose admissible monomials $M_{ve}$ for all the
adjacent edges which transform equivariantly with respect to the
action of $D(\Gamma)$.  Then $D(\Gamma)$ acts on the associated
equations of splice type so long as the higher order terms in these
equations transform appropriately under the action of the group.

\begin{definition}\label{splicequotient} Let $\Gamma$ be a \minimal{}
  resolution tree satisfying the semigroup and congruence conditions.
  Let $\Delta$ be the corresponding splice diagram, $z_{w}$ a variable
  associated to each leaf $w$, $M_{ve}$ an admissible monomial for each
  node $v$ and adjacent edge $e$ satisfying the
  $D(\Gamma)$--equivariance condition.  Then \emph {splice diagram
    equations for $\Gamma$} are equations of the form
$$\sum_{e}a_{vie}M_{ve} +H_{vi}=0,\quad i=1,\ldots,\delta_{v}-2,\quad v
\text{ a node}\,,$$
where
\begin{itemize}
\item for every $v$, all maximal minors of the matrix $(a_{vie})$ have
  full rank;
\item $H_{vi}$ is a convergent power series in the $z_{w}$'s all of
  whose monomials have $v$--weight $>d_{v}$;
\item for each $v$, the monomials in $H_{vi}$ transform under
  $D(\Gamma)$ in the same way as do the $M_{ve}$'s.
\end{itemize}
\end{definition}

We are ready for the careful statement of the main result of this
paper.

\begin{theorem}[Splice-quotient singularities]\label{th:main}
Suppose $\Gamma$ is \minimal{} and satisfies the semigroup and the
  congruence conditions. Then:
\begin{enumerate}
\item Splice diagram equations for $\Gamma$ define an isolated
complete intersection singularity $(X,o)$.
\item The discriminant group $D(\Gamma)$ acts freely on a
punctured neighborhood of $o$ in $X$.
\item $Y=X/D(\Gamma)$ has an isolated normal surface singularity,
and a good resolution whose associated dual graph is $\Gamma$.
\item $X\rightarrow Y$ is the universal abelian covering.
\item For any node $v$, the $v$--grading on functions on $Y$ (induced
  by the $v$--grading on $X$) is $\disc(\Gamma)$ times the grading by
  order of vanishing on the exceptional curve\comment{Addedc 70/02}
  $E_v$.\label{item:main-grading}
\item $X\rightarrow Y$ maps the curve $z_{w}=0$ to an irreducible
  curve, whose proper transform on the good resolution of $Y$ is
  smooth and intersects the exceptional curve transversally, along
  $E_{w}$. In fact the function $z_w^{\disc(\Gamma)}$, which is\comment{added the
  order of vanishing}
  $D(\Gamma)$--invariant and hence defined on $Y$, vanishes to order
  $\disc(\Gamma)$ on this curve.\label{item:main-endcurve}
\end{enumerate}
\end{theorem}

The first assertion of the theorem has already been proved in Theorem
\ref{th:splice is CI}.  For the second, recall (Proposition
%\ref{prop:free off hyperplanes})
\ref{prop:5.3}) that the fixed locus of a
non-identity element of $D(\Gamma)$ is contained in some subspace
$z_{w}=z_{w'}=0$, which intersects the germ of $X$ only at the origin
(Corollary \ref{isolated}).
It follows that $Y$ has an isolated normal singularity, and the main
point is to show the resolution dual graph equals $\Gamma$.  Once that
is achieved, the link of $Y$ will be a rational homology sphere whose
universal abelian covering has order $\disc(\Gamma)$, hence must be
given by the abelian covering provided by the link of $X$; this gives
the fourth assertion.

The fifth assertion follows immediately from the sixth if
one\comment{Added 07/02}
restricts to monomials (although its general proof will involve more
  work).
\begin{lemma}\label{le:grading}
  Let $d=\disc(\Gamma)$. Assuming assertion (\ref{item:main-endcurve})
  of the theorem, if $z_w$ is the variable corresponding to a leaf
  then the order of vanishing of $z_w^{d}$ on $E_v$ is
  $\ell_{vw}$ ($z_w^{d}$ is $D(\Gamma)$--invariant, so it
  is defined on $Y$).
\end{lemma}
\begin{proof}
  The divisor of the function $z_w^d$ on $\bar Y$ has the form
  $$(z_w^d)=dD+\sum r_u E_u\,,\quad\text{sum over vertices $u$ of
    $\Gamma$ other than }w\,,$$
  Where $D$ is the proper transform of
  the curve given by $z^d_w=0$ in $Y$.  This divisor dots to zero with
  each $E_{v}$, whence $r_v$ is the $vw$ entry of $-dA(\Gamma)^{-1}$.
  Theorem \ref{th:props} then says $r_v=\ell_{vw}$; the $v$--weight of
  $z_w$.
  \end{proof}

By taking a closer look at the weighted blow-up as in Lemma
\ref{le:tilde} we will first check the theorem in the one-node case of
\cite{neumann83} (giving a new proof for this case), and then proceed
by induction on the number of nodes. This induction involves choosing
an end-node $v^*$ of $\Gamma$ (or $\Delta$), and reducing to $\tilde
\Gamma$, obtained by removing $v^*$ and its adjacent strings leading
to leaves. But if there are no curves between $v^*$ and the node in
the remaining direction, then we must first blow-up between these
nodes, and create a new $-1$ curve.  This will guarantee that the new
$\tilde \Gamma$ is also \minimal{} (and that there is a leaf
corresponding to the removed $v^{*}$).  It is easy to check the
following:

\begin{lemma} Suppose $\Gamma$ satisfies the semigroup and congruence
conditions, and has adjacent
nodes.  Let $\Gamma'$ be the graph obtained from blowing-up between
the adjacent nodes.  Then $\Gamma$ and $\Gamma'$ have the same
splice diagram and representation of the discriminant group on the
space of ends.  In particular, splice diagram equations for
$\Gamma$ are the same for $\Gamma'$.\qed
 \end{lemma}
  We also need the \minimal{} version of a well known and
 classical\comment{WDN: Another minor change.} lemma:

\begin{lemma}\label{le:cyclic quotient}
  Consider a negative-definite string of rational curves
  $$\xymatrix@R=6pt@C=24pt@M=0pt@W=0pt@H=0pt{\\
    \overtag{\Circ}{-b_1}{8pt}\lineto[r]&
    \overtag{\Circ}{-b_2}{8pt}\dashto[r]&
    \dashto[r]&\overtag{\Circ}{-b_k}{8pt}\\&~}$$
  where either $b_{i} >1$ for
  $i=1,\dots,k$, or
    $k=1$ and $b_1=1$.  Write the
  continued fraction
    $$n/p=b_{1}-1/b_{2}-1/\cdots -1/b_{k}\,.$$
    Then the cyclic
    quotient singularity $\C ^{2}/\langle T\rangle$, where
    $$T(x,y)=\bigl(\exp (2\pi i/n) x,\,\exp (2\pi ip/n)y\bigr)\,,$$
    has a resolution with the above string of exceptional curves.  The
    proper transform of the image of $x=0$ on this resolution
    intersects once transversally on the right.\qed
\end{lemma}
Reading the string in the other direction yields $n/p'$, where
$pp'\equiv 1 \pmod{n}$.  Note $n=p=p'=1$ is allowed; in all other
cases, $n>1$ and $p<n$.

  \section{The case of one node}\label{sec:one node}

Recall that if $\C ^{t}$ has coordinates $z_{i}$ of weight
$m_{i}$, the weighted blow-up
$Z\rightarrow \C ^{t}$ has an open covering $U_{i}$,
each of
which is a quotient of an affine space $V_{i}$ by a cyclic group of order
$m_{i}$.  $V_{1}$ has coordinates $A_{1},\dots,A_{t}$,
 related to the $z_{i}$ via
$$z_{1}=A_{1}^{m_{1}},~ z_{2}=A_{1}^{m_{2}}A_{2},~ \dots~,~
z_{t}=A_{1}^{m_{t}}A_{t};$$ $U_{1}$ equals $V_{1}$
modulo the action of the
cyclic
group generated by $$S=[-1/m_{1},m_{2}/m_{1},\dots,m_{t}/m_{1}]\,.$$
A finite group $D$ of diagonal matrices acting on $\C ^{t}$ preserves the
weight filtration, so it lifts to an action
on $Z$, and one has a proper birational map $Z/D \rightarrow \C
^{t}/D$. A diagonal $[\beta_{1},\dots,\beta_{t}]$ acting
on $\C ^{t}$ may be lifted to one acting on $V_{1}$ via
$$-\beta_{1}S+[0,\beta_{2},\dots,\beta_{t}];$$
this lift depends on the choice of $\beta_{1}$.  The naturally
defined group
$D_{1}$, generated by $S$ and all lifts of elements of $D$, acts on
$V_{1}$, and induces an isomorphism $V_{1}/D_{1} \simeq U_{1}/D$ onto
an open subset of $Z/D$.  In
our situation, $D_{1}$ will contain a pseudo-reflection of the
form $[1/d,0,\dots,0]$. Dividing $V_{1}$ first by this action
produces another affine space $\bar{V}_{1}$, with coordinates
$\bar{A}_{1},A_{2},\dots,A_{t}$, where $\bar{A}_{1}=A_{1}^{d}$;
the quotient group $\bar{D}_{1}$ acts on $\bar{V}_{1}$ with
quotient $U_{1}/D$.   Note that the image of the exceptional divisor
is given by
$\bar{A}_{1}=0$.

Consider now a \minimal{} resolution graph $\Gamma$ with one node,
given by the diagram
$$\xymatrix@R=4pt@C=24pt@M=0pt@W=0pt@H=0pt{\\
\lefttag{\Circ}{n_2/p_2}{8pt}\dashto[ddrr]&
&\hbox to 0pt{\hss\lower 4pt\hbox{.}.\,\raise3pt\hbox{.}\hss}
&\hbox to 0pt{\hss\raise15pt
\hbox{.}\,\,\raise15.7pt\hbox{.}\,\,\raise15pt\hbox{.}\hss}
&\hbox to 0pt{\hss\raise 3pt\hbox{.}\,.\lower4pt\hbox{.}\hss}
&&\righttag{\Circ}{n_{t-1}/p_{t-1}}{8pt}\dashto[ddll]\\
\\
&&\Circ\lineto[dr]&&\Circ\lineto[dl]\\
\lefttag{\Circ}{n_1/p_1}{8pt}\dashto[rr]&&
\Circ\lineto[r]&\overtag{\Circ}{-b}{8pt}\lineto[r]&\Circ
\dashto[rr]&&\righttag{\Circ}{n_{t}/p_{t}}{8pt}\\&~\\&~}
$$
The strings of $\Gamma$ are described uniquely by the continued
fractions shown, starting from the node (by the \minimal{}ity
condition of Definition \ref{minimality}).  Set $N=n_{1}n_{2}\dots
n_{t}$, $N_{i}=N/n_{i}$; then
$d=N(b-\sum_{i=1}^{t}p_{i}/n_{i})=\disc(\Gamma)$ is the determinant.
The splice diagram produces leaf variables $z_{i}$ of weight $N_{i}$
for $1\leq i \leq t$.

The affine space
$V_{1}$ above has coordinates $A_{i}$, $1\leq i \leq t$ related to the
$z_{i}$ via
$$z_{1}=A_{1}^{N_{1}},~z_{i}=A_{1}^{N_{i}}A_{i},\quad 2\leq i \leq t\,.$$
The quotient $V_{1}\rightarrow U_{1}$ comes from dividing by
$$S=[-1/N_{1},N_{2}/N_{1},\dots,N_{t}/N_{1}]\,.$$
The discriminant group $D(\Gamma)$ acts on
$\C^{t}$ and the weighted blow-up $Z$.  $D(\Gamma)$ is generated by
\begin{align*}
&[N_{1}/(n_{1}d) + p'_{1}/n_{1},~N_{2}/(n_{1}d),~\dots~,~
N_{t}/(n_{1}d)]\,,\\
&[N_{1}/(n_{2}d),~N_{2}/(n_{2}d) +
p'_{2}/n_{2},~\dots~,~N_{t}/(n_{2}d)]\,,  \\
&\qquad\qquad\cdots\qquad\qquad\cdots\qquad\qquad\cdots\qquad\qquad.
\end{align*}
A simple calculation verifies the following lifts of these elements
to $V_{1}$:
\begin{align*}
  T_{1}&=[1/(n_{1}d)+p'_{1}/N,~-p'_{1}/n_{2},~\dots~,~-p'_{1}/n_{t}]\\
T_{2}&=[1/(n_{2}d),~p'_{2}/n_{2},~0,~\dots~,~0]\\
&\quad\cdots\qquad\qquad\cdots\\
T_{t}&=[1/(n_{t}d),~0,~\dots,~ 0,~p'_{t}/n_{t}]\,.
\end{align*}
Consider the lifted group $D_{1}$ generated by $S$
and these $T_{i}$. Writing $p_{1}p'_{1}=kn_{1}+1$, we may
replace $S$ and $T_{1}$ by
\begin{align*}
\bar{S}&=S^{p'_{1}}T_{1}^{n_{1}}=[1/d,~0,~\dots~,~0]\\
\bar{T}&=\,S^{k}T_{1}^{p_{1}}\,=[1/N+
p_{1}/(n_{1}d),~-1/n_{2},~\dots~,~-1/n_{t}]\,,
\end{align*}
and then replace $\bar{T}$ by
$$\tilde{T}=\bar{T}T_{2}^{p_{2}}\dots T_{t}^{p_{t}}=\bigl[1/N
+\sum_{i=1}^{t}p_{i}/(n_{i}d),~0,~\dots~,~0\bigr]= [b/d,~0,~\dots~,~0]\,.$$
Since $\tilde T$ is a power of $\bar S$,
$D_{1}$ is generated by the pseudoreflection
$\bar{S}$ and the $T_{i},i>1$.  Dividing $V_{1}$ by $\bar{S}$ gives the affine
space $\bar{V}_{1}$, with coordinates
$\bar{A}_{1}=A_{1}^{d},A_{2},\dots,A_{t}$, and group action
generated by
\begin{align*}
  \bar{T}_{2}&=[1/n_{2},~p'_{2}/n_{2},~0,~\dots~,~0]\\
&\qquad\qquad\cdots\\
\bar{T}_{t}&=[1/n_{t},~0,~\dots~,~0,~p'_{t}/n_{t}]\,.
\end{align*}
 In this one-node case, the admissible
 monomials are
$z_{i}^{n_{i}}$, $i=1,\dots,t$,
and the splice diagram equations are of the form
$$\sum_{j=1}^{t}a_{ij}z_{j}^{n_{j}} +\ H_{i}(\underline z)=0,\
i=1,\dots, t-2\,.$$
A monomial $\prod z_{j}^{\alpha_{j}}$ is allowed to appear in
one of the
convergent power series $H_{i}$ iff it transforms under
$D(\Gamma)$ as do the admissible monomials, and if the $v$--weight $\sum
\alpha_{j}N_{j}$ is $>N$.  By the assumption on the coefficient
matrix $(a_{ij})$, one may take linear combinations of the series
and rewrite as
$$z_{i}^{n_{i}}+a_{i}z_{t-1}^{n_{t-1}}+b_{i}z_{t}^{n_{t}}+
H_{i}(\underline z)=0,\ \ 1\leq i\leq t-2$$
with $a_ib_j-a_jb_i\ne0$ for all $i\ne j$, and all $a_i$ and $b_i$
nonzero. The associated graded ring
with respect to the weight filtration is the familiar
Brieskorn complete intersection (except that some $n_{i}=1$ is
possible).

As in Section \ref{sec:blow-up}, one finds the proper transform
$X_{1}$ on $V_{1}$ defined by:
\begin{align*}
  1+a_{1}A_{t-1}^{n_{t-1}}+b_{1}A_{t}^{n_{t}}+
A_{1}J_{1}(\underline A)&=0\\
A_{i}^{n_{i}}+a_{i}A_{t-1}^{n_{t-1}}+b_{i}A_{t}^{n_{t}}+
A_{1}J_{i}(\underline A)&=0,\quad 2\leq\ i\leq t-2\,.
\end{align*}
Here one has used the weight condition on the $H_{i}$ to define
$$H_{i}(A_{1}^{N_{1}},A_{1}^{N_{2}}A_{2}
,\dots,A_{1}^{N_{t}}A_{t})=A_{1}^{N+1}J_{i}(\underline A)\,.$$
Intersecting with the exceptional divisor $A_{1}=0$ gives the
polynomials
equations
\begin{equation}\label{eq:blowup1}
  \begin{split}
1+a_{1}A_{t-1}^{n_{t-1}}+b_{1}A_{t}^{n_{t}}&=0\,,\\
A_{i}^{n_{i}}+a_{i}A_{t-1}^{n_{t-1}}+b_{i}A_{t}^{n_{t}}
&=0\,,\quad 2\leq\ i\leq t-2\,.
\end{split}
\end{equation}
This is a smooth curve, so $X_{1}$ is smooth.  Since $\bar{S}$
leaves $A_{2},\dots,A_{t}$ invariant and acts equivariantly on
the equations, it follows that each term
$A_{1}J_{i}(A_{k})$ must be invariant, and hence is a
power series
$\bar{A}_{1}\bar{J}_{i}$
in $\bar{A}_{1},A_{2},\dots,A_{t}$.  So the quotient
$\bar{X_{1}}$ of $X_{1}$ by $\bar{S}$ is defined on
$\bar{V}_{1}$ via
\begin{equation}\label{eq:blowup}
  \begin{split}
      1+a_{1}A_{t-1}^{n_{t-1}}+b_{1}A_{t}^{n_{t}}+
\bar{A}_{1}\bar{J}_{1}&=0\,,\\
A_{i}^{n_{i}}+a_{i}A_{t-1}^{n_{t-1}}+b_{i}A_{t}^{n_{t}}+
\bar{A}_{1}\bar{J}_{i}&=0\,,\quad 2\leq\ i\leq t-2\,.
  \end{split}
\end{equation}
We divide $\bar{X_{1}}$ by the action of $\bar{D}_{1}:=
D_{1}/\bar{S}$, which is generated by the images of $\bar{T}_{i},\
i=2,\dots ,t$.
First, the group acts transitively on the connected
components of $\bar{A}_{1}=0$, since, eg, every component contains a
point with $A_{t}=0$ (cf. Lemma \ref{le:4.1}); since $A_{1}=0$ in $X_{1}$
is smooth, the image of the exceptional divisor $\bar{A}_{1}=0$ is irreducible.
Next, the action is free off
$\bar{A}_{1}=0$, and fixed points occur exactly when another
coordinate is 0.  At a point of $\bar{X_{1}}$ where
$\bar{A}_{1}=A_{2}=0$, the above equations (\ref{eq:blowup}) determine
$A_{k}^{n_{k}}$ for $k=3,\dots,t$ uniquely. Thus there are
$n_{3}\dots n_{t}$ such points, and they are permuted by the subgroup
generated by $\bar{T}_{3},\dots,\bar{T}_{t}$.  So, there is one orbit
of such fixed points, and the stabilizer is generated by
$\bar{T}_{2}$.  The equations imply that $\bar{A}_{1},A_{2}$ are local
coordinates at such a point.

  We have now a familiar local picture as described in Lemma
  \ref{le:cyclic quotient} and following comments: divide $\C^{2}$ by
  the action $[1/n_{2},p'_{2}/n_{2}]$, resolve the cyclic quotient
  singularity according to the corresponding string in $\Gamma$ (if
  $n_2=1$ we blow-up once), and
  consider the proper transforms of the images of the two coordinates
  axes on $\C^{2}$.  We get a string of curves
$$\xymatrix@R=18pt@C=18pt@M=0pt@W=0pt@H=0pt{
&\lineto[dd]^(.8)C&\lineto[ddrr]&\lineto[dlll]&&&&\lineto[dd]^(.2){P}&\\
&&&&\dots&\dots&\dots&&\lineto[dll]^(.8){E_1}
\\
&&&&&&&&}
$$
The left-hand curve $C$ corresponds to the transform of
$\bar{A}_{1}=0$ (which will thus be a central exceptional curve in a
resolution of $X$), the right curve $P$ is the proper transform of
$A_{2}=0$, and the remaining curves form the string of exceptional
curves that resolve the cyclic quotient singularity.  The continued
fraction expansion from left to right is $n_{2}/p_{2}$.
$A_{2}^{n_{2}}$ vanishes $n_{2}$ times along $P$ and $p'_{2}$ times
along the adjacent exceptional curve $E_{1}$; and,
$\bar{A}_{1}^{n_{2}}$ vanishes once along $E_{1}$.

Putting this all together for the other $U_{i}$ and $V_{i}$, we conclude that
the quotient variety $X/D(\Gamma)$ has a resolution consisting of
a smooth central curve, and $t$ rational strings emanating from it
corresponding to $n_{i}/p_{i},\ 1\leq i\leq t$.  It remains to
show the central curve is rational, and its self-intersection is
$-b$ (equivalently, the determinant $\bar{d}$ of the intersection pairing
equals $d$).

Restricting to the exceptional divisor $\bar{A}_{1}=0$ in
$\bar X_1/\bar{D}_{1}=X_1/D_1$
gives new variables $C_{i}=A_{i}^{n_{i}},\ 2\leq i\leq t$, in which
the defining equations in (\ref{eq:blowup1}) become linear.  So, the
quotient is a line in the coordinate space, hence is rational.

Finally, note that $z_{2}^{dn_{2}}$ is invariant under the
discriminant group, hence is a function on $X/D(\Gamma)$; we consider
its proper transform in the minimal resolution of
$\bar{V}_{1}/\bar{D}_{1}$.  (The proper transform misses $V_{2}$, so
is completely contained in $V_{1}$.)  Note first that
$$z_{2}^{dn_{2}}=A_{1}^{dN}A_{2}^{dn_{2}}=\bar{A}_{1}^{N}A_{2}^{dn_{2}}\,.$$
By an earlier remark, $z_{2}^{dn_{2}}$ vanishes $N/n_{2}
+dp'_{2}$ times along $E_{1}$, and $dn_{2}$ times along $P$.
Thus the divisor of the function $z_2^{dn_2}$ has the form
\begin{equation}
  \label{eq:15}
(z_2^{dn_2})=dn_{2}P+\bigl(N/n_2
+dp'_{2}\bigr)E_{1}+\sum_{i>1}r_{i}E_{i}
\end{equation}
the last sum being over all the other exceptional curves in the
\minimal{} %good
resolution.  Since the divisor of a function dots to zero
with each $E_i$,  we see that
$$-\frac{1}{dn_2}\Bigl(\bigl(N/n_{2}
+dp'_{2}\bigr)E_{1}+\sum_{i>1}r_{i}E_{i}\Bigr)$$
represents the element
$e_1$ in the dual $\E^*$ of $\E = \bigoplus \Z \cdot E_{v}$. On the
other hand, $e_1=\sum_i(e_1\cdot e_i)E_i$, so
$$e_1\cdot e_1=-\bigl(N/n_2+dp'_{2}\bigr)/(dn_{2})\,.$$
Comparing this with the
value $e_1\cdot e_1=-N/(n_{2}^{2}\bar{d}) -p'_{2}/n_{2}$ of
Proposition \ref{prop:endwt} yields that $\bar{d}=d$, as desired.

So, the constructed singularity has a resolution with dual graph our
original $\Gamma$.  Since we have constructed an abelian covering
of degree $d$ equal to the discriminant of our singularity, the map
$X\rightarrow X/D(\Gamma)$ must be the universal abelian covering.

To prove the assertion (\ref{item:main-endcurve}) of the main
theorem, note that the proper
transform of $z_{i}^{n_{i}d}$ on the given \minimal{} %good
resolution
intersects the $i$-th quotient string, once transversally on the end,
with multiplicity $n_{i}d$.

Finally\comment{added 07/02}, for assertion (\ref{item:main-grading})
of the main theorem, the central exceptional curve in the resolution
of $X$ maps to the central curve $E_v$ in the resolution of $Y$, so
the gradings on functions on $Y$ (ie, $D(\Gamma)$--invariant
functions on $X$) given by order of vanishing on these curves agree up
to a constant multiple. Since the central curve of the resolution of
$X$ is the curve obtained by blowing up the $v$--grading, order of
vanishing on it is given by the $v$--grading. This proves
(\ref{item:main-grading}) up to a multiple; that the multiple
is correct is confirmed in Lemma \ref{le:grading}.

\section{The inductive procedure}\label{sec:induct}

Assume we have a set of splice diagram equations for a \minimal{} $\Gamma$, which
is assumed to have more than one node.  The inductive assumption is
that the theorem is true for a graph with fewer nodes.  As mentioned,
we may assume (after blowing-up) that any two nodes of $\Gamma$ have at least one curve
between them.  Let $\Delta$ be the splice
diagram and $v^{*}$ an end-node, viewed on the left side of the
diagrams:
\begin{gather*}
\xymatrix@R=4pt@C=24pt@M=0pt@W=0pt@H=0pt{\\
&\lefttag{\Circ}{n_\sigma/p_\sigma}{8pt}\dashto[ddrr]&
&&&&\\&&&&&&&&&& \\&
&&\Circ\lineto[dr]&&\hbox to 0 pt{$\scriptstyle n/p~\rightarrow$\hss}&&&&\Circ\lineto[dl]\dashto[ur]\\
\Gamma=\quad&
&\hbox to 0pt{\hglue-8pt\vbox to 0 pt{\vss\vss\vdots\vss}\hss}&&
\undertag{\overtag{\Circ}{-b}{8pt}}{v^*}{6pt}\lineto[r]&
\dashto[rr]&&\lineto[r]&\Circ&&
\hbox to 0 pt{\hss\vbox to 0 pt{\vss\vss\vdots\vss}\hglue5pt}
\\&&&\Circ\lineto[ur]&&&&&&\Circ\lineto[ul]\dashto[dr]\\
&&&&&&&&&&
\\&\lefttag{\Circ}{n_1/p_1}{8pt}\dashto[uurr]
}
\\
\xymatrix@R=6pt@C=30pt@M=0pt@W=0pt@H=0pt{\\
&\Circ\lineto[ddrr]^(.7){n_\sigma}&&&&&&&\\&&&&&&&\dashto[ur]
\\
\Delta=&&\hbox to 0pt{\hglue-8pt\vbox to 0 pt{\vss\vss\vdots\vss}\hss}
&\undertag\Circ{v^*}{6pt}\lineto[rrr]^(.2)r^(.8)s&&&
\Circ\lineto[ur]^(.6){m_\mu}\lineto[dr]_(.6){m_1}
&\hbox to 0 pt{\hss\vbox to 0 pt{\vss\vss\vdots\vss}\hglue-8pt}&
\\&&&&&&&\dashto[dr]\\&\Circ\lineto[uurr]_(.7){n_1}&&&&&&&
\\&}
\end{gather*}
The continued fractions from $v^{*}$ to the leaves are given by
$n_{i}/p_{i},\ 1\leq i\leq \sigma$, and from $v^{*}$ to the adjacent
node by $n/p$.  Let $N=n_{1}\dots n_{\sigma},
N_{i}=N/n_{i}, M=m_{1}\dots m_{\mu}$.  One has the determinant
calculation
\begin{equation}
  \label{eq:6}
s=Nn(b-\sum_{i=1}^{\sigma}p_{i}/n_{i}-p/n)\,.
\end{equation}
Moreover, if $d=\disc(\Gamma)=\det(-A(\Gamma))$,
we have the
relation (see Proposition \ref{prop:edge det})
\begin{equation}
  \label{eq:7}
rs-MN=dn\,.
\end{equation}
Induction will involve $\tilde\Gamma $, obtained
by removing from $\Gamma$ the vertex corresponding to $v^{*}$ plus
the $\sigma$ strings of rational curves to the leaves.
The corresponding splice
diagram $\tilde\Delta$ has a new leaf $w^{*}$ in the location of
$v^{*}$, but one loses the
$\sigma$ leaves of $\Delta$ adjacent to $v^{*}$.  The left sides of
$\tilde\Gamma$ and $\tilde\Delta$ are
\begin{gather*}
\xymatrix@R=4pt@C=24pt@M=0pt@W=0pt@H=0pt{\\
&&&&&& \\&&\hbox to 0 pt{\hss$\scriptstyle n/p~\rightarrow$\hss}
&&&\Circ\lineto[dl]\dashto[ur]\\
\tilde\Gamma=\quad&
\Circ\lineto[r]&
\dashto[r]&
\lineto[r]&\Circ
&&\hbox to 0
pt{\hss\vbox to 0 pt{\vss\vss\vdots\vss}\hglue5pt}
\\&&&&&\Circ\lineto[ul]\dashto[dr]\\
&&&&&&
}
\\
\xymatrix@R=6pt@C=30pt@M=0pt@W=0pt@H=0pt{\\
&&&&&\\&&&&\dashto[ur]
\\
\quad\tilde\Delta=
&\undertag{\Circ}{w^*}{4pt}\lineto[rr]^(.75)n&&
\Circ\lineto[ur]^(.6){m_\mu}\lineto[dr]_(.6){m_1}
&\hbox to 0 pt{\hss\vbox to 0 pt{\vss\vss\vdots\vss}\hglue-8pt}&
\\&&&&\dashto[dr]\\&&&&&
\\&}
\end{gather*}
For edges of $\tilde\Delta$ which
point away from $w^{*}$, the weights are the same as they were for
$\Delta$. For the edge pointing toward $w^*$ at a node $v$ the new
edge weight is now given by (see Lemma \ref{le:10.7}):
\begin{equation}
  \label{eq:10}
  \tilde{d}_{v1}=\frac1d\left(rd_{v1}-N(d_{v}/d_{v1})
  (\ell'_{w^{*}v})^{2}\right)\,.
\end{equation}
Note also that
\begin{equation}
  \label{eq:12}
  r=\disc(\tilde\Gamma)=\det(-A(\Gamma))\,.
\end{equation}
We use exactly the same notation as in Section \ref{sec:blow-up} for the
variables and equations for the singularity $X(\Delta)$, including the
$v^{*}$--blow-up to $X_{1}\subset V_{1}$.  The coordinates
corresponding to the left leaves are $x_1,\dots,x_\sigma$ and those
corresponding to the right leaves are $y_1,\dots, y_\tau$. The
coordinate space $V_1$ has coordinates $A_i,\z_j$, related to the
$x_i,y_j$ by
\begin{equation}
  \label{eq:5}
x_1=A_1^{N_1r};~x_i=A_1^{N_ir}A_i,\quad i=2,\dots,\sigma;\qquad
y_j=A_1^{N\ell'_j}\z_j,\quad j=1,\dots,\tau\,,
\end{equation}
where we are abbreviating $\ell'_j:=\ell'_{v^*w_{\sigma+j}}$ (the reduced
linking number of $v^{*}$ and one of the $\tau$ outer leaves;
similarly, we will write $\ell_{jk}$ for the linking numbers between
two such leaves).  On $V_{1}$, the cyclic group action is generated by
$$S=[-1/rN_{1},n_{1}/n_{2},\dots,n_{1}/n_{\sigma};
\ n_{1}\ell'_{1}/r,\dots]\,.$$
We lift the discriminant group $D(\Gamma)$ to a group $D_{1}$ (of
order $dN_1r$) acting on $V_{1}$ and, by equivariance of the action, on
$X_{1}$.

The leaf-generators of $D(\Gamma)$ acting in $x_i,y_j$ coordinates are
written as follows, where a semi-colon distinguishes the first
$\sigma$ entries from the last $\tau$:
\begin{align*}
&[rN_{1}/(n_{1}d)+p'_{1}/n_{1},~rN_{2}/(n_{1}d),~\dots~,
rN_{\sigma}/(n_{1}d);~\ell'_{1}N_{1}/d,~\dots]\\
&[rN_{1}/(n_{2}d),~rN_{2}/(n_{2}d)+p'_{2}/n_{2},~\dots~,
rN_{\sigma}/(n_{2}d);~\ell'_{1}N_{2}/d,~\dots]\\
&\qquad\qquad\qquad\cdots\quad\cdots\\
&[rN_{1}/(n_{\sigma}d),~rN_{2}/(n_{\sigma}d),~\dots~
,~rN_{\sigma}/(n_{\sigma}d)+p'_{\sigma}/n_{\sigma};~\ell'_{1}N_{\sigma}/d,~\dots]\\
&[\ell'_{1}N_{1}/d,~\dots~,~\ell'_{1}N_{\sigma}/d;~
\ell_{11}/d,~\ell_{12}/d,~\dots]\\
&\qquad\qquad\qquad\cdots\quad\cdots\\
&[\ell'_{\tau}N_{1}/d,~\dots~,~\ell'_{\tau}N_{\sigma}/d;~
\ell_{\tau1}/d,~\ell_{\tau2}/d,~\dots]
\end{align*}
One verifies by substituting from (\ref{eq:5}) that the
following give lifts of these generators
of $D(\Gamma)$ to the coordinates of $V_{1}$:
\begin{align*}
T_{1}&=[1/(n_{1}d)+p'_{1}/(rN),~-p'_{1}/n_{2},~\dots~,~-p'_{1}/n_{\sigma};~
-p'_{1}\ell'_{1}/r,~\dots]\\
T_{2}&=[1/(n_{2}d),~p'_{2}/n_{2},~0,~\dots~,~0;~0,~\dots ]\\
&\qquad\qquad\qquad\cdots\quad\cdots\\
T_{\sigma}&=[1/(n_{\sigma}d),~0,~\dots~,~0,~p'_{\sigma}/n_{\sigma};~
0,~\dots ]\\
R_{1}&=[\ell'_{1}/(rd),~0,~\dots~,~0;~\ell_{11}/d
-N(\ell'_{1})^{2}/(rd),~\ell_{12}/d-N\ell'_{1}\ell'_{2}/(rd),~\dots ]\\
&\qquad\qquad\qquad\cdots\quad\cdots\\
R_{\tau}&=[\ell'_{\tau}/(rd),~0,~\dots~,~0;~\ell_{\tau1}/d
-N\ell'_{\tau}\ell'_1/(rd),~\ell_{\tau2}/d-N\ell'_{\tau}\ell'_{2}/(rd),~\dots ]
\end{align*}
The lifted discriminant group $D_{1}$ is generated by $S$ and the
$T_{i}$ and $R_{j}$.  Writing $p_{1}p'_{1}=kn_{1}+1$, replace the
generators $S$ and $T_{1}$ by new generators
\begin{align*}
  S'&=S^{p'_{1}}T_{1}^{n_{1}}=[1/d,~0,~\dots~,~0;~0,~\dots ]\\
T'&=S^{k}T_{1}^{p_{1}}=[1/(rN)+p_{1}/(n_{1}d),~
-1/n_{2},~\dots~,~-1/n_{\sigma};~ -\ell'_{1}/r,~\dots ]\,.
\end{align*}
Also replace $T'$ by $\tilde{T}$, where
$$\tilde{T} ^{-1}=T'T_{2}^{p_{2}}\dots
T_{\sigma}^{p_{\sigma}}=[1/(rN)
+\sum_{i=1}^{\sigma}p_{i}/(dn_{i}),~0,~\dots~,~0;~-\ell'_{1}/r,~\dots ]\,.$$
Recall (Lemma \ref{le:4.1})
that $X_{1}$ is smooth save for the $N_{1}$ points on the exceptional curve
$A_{1}=0$ where all $\z_{k}$ are 0; note that the group generated by
$T_{2},\dots,T_{\sigma}$ acts transitively on them, so dividing by
the action will give a connected exceptional curve.
$D_{1}$ also acts on $X_{1}$.

Divide $V_{1}$ by the pseudo-reflection $S'$, giving a degree $d$
covering $V_{1}\rightarrow \bar{V}_{1}$, which is unramified off the
divisor $A_{1}=0$.  In the new affine space $\bar{V}_{1}$, we have
coordinates $\bar{A}_{1}=A_{1}^{d},\
A_{2},\dots,A_{\sigma},\z_{1},\dots,\z_{\tau}$.

The action of the quotient $\bar{D}_{1}=D_{1}/\langle S'\rangle$ on
$\bar V_1$ is generated by multiplying the first entries of the old
generators by $d$, yielding new generators
\begin{align*}
  \bar{T}&=[-d/(rN)-\sum_{i=1}^{\sigma}
p_{i}/n_{i},~0,~\dots~,~0;~\ell'_{1}/r,~\dots ]\\
\bar{T}_{2}&=[1/n_{2},~ p'_{2}/n_{2},~0,~\dots~,~0;~ 0,~\dots ]\\
&\qquad\qquad\qquad\cdots\quad\cdots\\
\bar{T}_{\sigma}&=[1/n_{\sigma},~0,~\dots~,~p'_{\sigma}/n_{\sigma};~
0,~\dots ]\\
\bar{R}_{1}&=[\ell'_{1}/r,~ 0,~\dots~,0;~ \ell_{11}/d
-N(\ell'_{1})^{2}/(rd),~\ell_{12}/d -N\ell'_{1}\ell'_{2}/(rd),~\dots ]\\
&\qquad\qquad\qquad\cdots\quad\cdots
\end{align*}
Equations (\ref{eq:6}) and (\ref{eq:7}) imply
$$-d/(rN) -\sum _{i=1}^{\sigma} p_{i}/n_{i}=M/(rn) +p/n\ -b\,,$$
whence the generator $\bar{T}$ above may be rewritten as
$$\bar{T}=[M/(rn) +p/n,~ 0,~\dots~,~0;~ \ell'_{1}/r,~
\ell'_{2}/r,~\dots ]\,.$$
The quotient $\bar{X_{1}}$ of $X_{1}$ by $S'$ is defined by the same
equations as $X_{1}$ (see equations (\ref{eq:transformed general splice})),
except that $A_{1}^{d}$ is replaced by $\bar{A_{1}}$. $\bar{X_{1}}$
and the exceptional divisor are smooth except at the $N_{1}$ points with all
$\z_{k}=0$, and they form one $\bar{D_{1}}$--orbit.  As in Section
\ref{sec:blow-up}, at such a point, the coordinates $A_{i}$ for $i>1$
are non-zero, and can be solved locally as convergent power series in
$\bar{A_{1}},\z_{1},\dots ,\z_{\tau}$ with non-zero constant term.
Thus, around each one of these singular points $\bar{X_{1}}$ is
described by equations in those variables.

\begin{lemma} These equations are splice diagram equations for $\tilde\Delta$.
\end{lemma}
\begin{proof}  The proof is essentially the same as Lemma
  \ref{le:tilde}: we verify that the equations have
  appropriate admissible monomials at every node, and all the other
  terms have higher order. The edge-weights of $\tilde\Delta$ are the
  same as those for $\Delta$ in directions away from $v^{*}$, while at
  a node $v$, the new edge weight $\tilde d_{vi}$ is $1/d$ times the
  value it had in Section \ref{sec:blow-up} (compare equations
  (\ref{eq:11}) and (\ref{eq:10})).  Monomials involve powers of
  $\bar{A}_{1}$ rather than $A_{1}$, so the ``$Q$'' term in the proof
  of Lemma \ref{le:tilde} is also divided by $d$. The right side
  of equation (\ref{eq:4}) is therefore now also divided by $d$;
  otherwise the proof goes through almost verbatim.
\end{proof}

Still considering one such singular point of $\bar{X_{1}}$, the stabilizer
is the subgroup
$G$ of $\bar{D}_{1}$ generated by
$\bar{T},\bar{R}_{1},\dots,\bar{R}_{\tau}$.\comment{we needed bars
on R's, right?}
This group acts effectively on the
space whose coordinates are $\bar{A}_{1}, \z_{1},\dots ,\z_{\tau}$.
Rewrite the group generators
as $\tau +1$--tuples:
\begin{align*}
  &[M/(rn) +p/n,~  \ell'_{1}/r,~ \ell'_{2}/r,~\dots ]\,,\\
&[\ell'_{i}/r,~\ell_{i1}/d-N\ell'_{1}\ell'_{i}/(rd),~
\ell_{i2}/d -N\ell'_{2}\ell'_{i}/(rd),~\dots ]\,,\quad i=1,\dots,\tau\,.
\end{align*}

\begin{lemma} The $\tau +1$ generators above are exactly the
leaf generators for the discriminant group associated to the
resolution diagram $\tilde\Gamma$.  In particular, $G$ maps
isomorphically onto $D(\tilde\Gamma)$, viewed as a diagonal subgroup
of $(\C^*)^{\tau +1}$.
\end{lemma}
\begin{proof} The only ``new'' weights
  on edges of $\tilde\Delta$ are those pointing towards $w^{*}$; other
  edge-weights are the same for $\Delta$ and $\tilde\Delta$.  Recall
  that $r=\disc(\tilde\Gamma)$. By Proposition \ref{prop:endwt} and
  Corollary \ref{cor:action}, the first generator above is exactly the
  leaf generator of $D(\tilde\Gamma)$ corresponding to $w^{*}$.
  Similarly, the first entries of the remaining $\tau$ generators
  above are exactly those of the corresponding leaf generators for
  $D(\tilde\Gamma)$.

  So, we need to compute the last $\tau$ entries for the last $\tau$
  generators.  Denote the linking numbers for $\tilde\Delta$ between
  pairs of these as $\tilde{\ell}_{ij}$ (keeping in mind the special
  definition when $i=j$ as in Proposition \ref{prop:endwt}).  To prove
  the lemma, we must show that for all $i$ and $j$
  $$\tilde{\ell}_{ij}/r=\ \ell_{ij}/d\ -N\ell'_{i}\ell'_{j}/rd\,.$$
  When $i=j$ this equation is just a special case of Lemma
  \ref{le:10.7}, so assume $i\ne j$. Denote by $v$ the vertex where the
  paths from $w^*$ to $w_i$ and $w_j$ diverge. Then
  $$\tilde\ell_{ij}=\frac{\ell_{ij}}{d_{v1}}\tilde d_{v1}=
  \frac{\ell_{ij}}{dd_{v1}}\left(rd_{v1}-N(d_v/d_{v1})
    (\ell'_{w^*v})^2\right)\,.$$
  Substituting this in the equation to
  be proved reduces it to the equation
  $\ell'_i\ell'_j=\frac{\ell_{ij}}{d_{v1}}(d_v/d_{v1})(\ell'_{w^*v})^2$.
  This is an easily checked equation between products of splice
  diagram weights.
\end{proof}

The last two lemmas together imply that
$Y_{1}:= \bar{X}_{1}/\bar{D}_{1}$ has a singular point which is
a quotient
of splice diagram equations for $\tilde\Delta$ by
$D(\tilde\Gamma)$.  By induction, this singularity has a resolution dual
graph $\tilde\Gamma$, and one knows where the proper transform of the
exceptional fiber intersects the diagram (ie, in the location of
$w^{*}$).  The only other singular points of
$Y_{1}$ come from fixed points of the group action.

\begin{lemma} $\bar{D}_{1}$ acts transitively on the points of
  $\bar{X}_{1}$ with $\bar{A_{1}}=A_{2}=0$.  At such a point,
  $\bar{A_{1}}$ and $A_{2}$ are local analytic coordinates, the
  stabilizer is generated by $\bar{T}_{2}$, and the image of the point
  on $Y_{1}$ is an $n_{2}/p_{2}$ cyclic quotient singularity.
\end{lemma}
\begin{proof} As before, the equations (\ref{eq:transformed general splice})
  (in which now $A_1$ only occurs to $d$-th powers which have been
  replaced by $\bar A_1$) show that no other coordinates can be
  0, and that $\bar{A_{1}}$ and $A_{2}$ are local coordinates (see
  also Lemma \ref{le:4.1}).

  Recall that the order of $\bar D_{1}$ is $d rN_{1}/d\ =rN_{1}$ and
  that the order of $G$ is $\det(-A(\tilde\Gamma))=r$.  The subgroup $G$ of
  $\bar D_1$, viewed as acting on the coordinates $\bar
  A_{1},\z_{1},\dots,\z_{\tau}$, contains no pseudoreflections (it is
  a discriminant group, so apply Proposition \ref{prop:5.3});
  therefore, no element of $G$ stabilizes the point
  in question.  Further, the subgroup generated by
  $\bar{T_{3}},\dots,\bar{T_{\sigma}}$, which has order $n_{3}\dots
  n_{\sigma}=N/n_{1}n_{2}$, acts on the points in question but only
  changing their entries in the slots $A_{3},\dots,A_{\sigma}$.  In
  particular, $G$ and the $\bar{T_{i}}$ for $i>2$ generate a subgroup
  of index $n_{2}$ in $\bar D_1$ which acts freely on the points, and
  the stabilizer of the point is generated by $\bar{T_{2}}$. Thus
  there are $rn_3\dots n_\sigma$ points in the orbit, and their image
  is an $n_2/p_2$ cyclic quotient singularity.  Finally, to see that
  $\bar D_1$ acts transitively on the points in question, we must show
  there are $rn_{3}\dots n_{\sigma}$ points with
  $\bar{A_{1}}=A_{2}=0$; but this is given in Lemma \ref{le:4.1}.
\end{proof}

Now, the exceptional divisor of $Y_{1}$ is connected, and the last two
lemmas show it is analytically irreducible at the singular points;
thus, the exceptional divisor is itself irreducible.  Take the
resolution of the above $n_{2}/p_{2}$ cyclic quotient singularity on
$Y_1$ (if $n_{2}=1$, this means blow-up a smooth point); then, as in
section \ref{sec:one node}, we get a string of rational curves
starting from the curve $\bar A_1=0$ with continued fraction
$n_2/p_2$, and the proper transform of the curve $A_2=0$ is smooth
and intersects
transversally in one point the end curve of this string.  Since the
curve $A_2=0$ arises from the curve $x_2=0$, this agrees with the last
part (\ref{item:main-endcurve}) of Theorem \ref{th:main}.

The action of $\bar{D_{1}}$ on the curve $\bar{A_{1}}=0$ reduces to
the direct product of $G$ acting in the $B$ coordinates and the group
generated by $\bar T_2,\dots,\bar T_\sigma$ acting in the $A$
coordinates. It thus acts freely except on the orbit of $\tilde
\Delta$ splice diagram singularities (which occur at the points where
all the $B_k$ are zero), and except on points where some $A_{i}=0$
($i>1$). The latter lead to cyclic quotient singularities on dividing
by $\bar D_1$.  The quotient $Y_1=\bar X_1/\bar D_1$ thus has one
splice-quotient singularity and, for each $i=2,\dots,\sigma$, one
$n_i/p_i$ cyclic quotient singularity, and is otherwise smooth along
its (irreducible) exceptional divisor.  We have already verified that
the image of $x_{i}=0$ (for $i>1$) on $Y$ has proper transform on
$Y_{1}$ vanishing correctly on the desired end-curves.

One can see the rest of the $v^{*}$ blow-up of $X$ by inverting
$x_{2}$, since the origin is an isolated point in the locus
$x_{1}=x_{2}=0$ (Corollary \ref{isolated}).  The resulting space
$Y_{2}$ adds an $n_{1}/p_{1}$ cyclic quotient (but misses the
$n_{2}/p_{2}$ one).  Denote by $\bar{Y}$ the union of $Y_{1}$ and
$Y_{2}$.  There is a partial resolution $\bar{Y}\rightarrow
X/D(\Gamma)$ whose exceptional curve is irreducible, and along which
sit (in a known way) $\sigma$ cyclic quotient singularities and a
splice-quotient corresponding to $\tilde\Gamma$.  Taking the
resolutions of all the quotient singularities as well as the
splice-quotient (which by the induction assumption has resolution
graph $\tilde\Gamma$) gives a resolution
$\tilde{Y}\rightarrow\bar{Y}\to X/D(\Gamma)$, with resolution dual
graph \emph{almost} guaranteed to be exactly $\Gamma$. The only point
to check is that the proper transform of the exceptional curve of
$\bar{Y}$ in $\tilde{Y}$ is a smooth rational curve of the correct
self-intersection.  But that is achieved simply by repeating the
entire procedure at a different end-node $v'$ of $\Gamma$; the curve
in question is then seen as a part of the resolution dual graph of
type $\Gamma'$ ($\Gamma$ minus $v'$ and its adjacent strings), hence
has the desired properties.

The assertion about the proper transform on $\bar{Y}$ of the image of
$y_{j}=0$ follows by induction, by considering the role of $B_{j}=0$
for the splice-quotient corresponding to $\tilde\Gamma$.

\smallskip We have proved all but part (\ref{item:main-grading}) of
Theorem \ref{th:main}. For the node $v^*$ that we just blew up, this
part follows by the same argument as in Section \ref{sec:one node} for
the one-node case.  For any other node $v$, we proceed by induction,
comparing the factorizations $X_{1}\rightarrow X \rightarrow Y$ and
$X_{1} \rightarrow \bar X_{1} \rightarrow Y_{1} \rightarrow Y$, and
the relevant valuations and weightings at various points.  Recall
(Lemma \ref{le:tilde}) that $X_{1}$ has $N_{1}$ splice-type singularities for
the splice diagram $\Delta '$, obtained by deleting from $\Delta$ the
node $v^{*}$ and its outer leaves, replacing by a new leaf $w^{*}$,
and adjusting weights as in equation (\ref{eq:11}).  (We had
previously denoted this diagram by $\tilde \Delta$; but in the current
proof, that notation is already being used, and means that the weights
in (\ref{eq:11}) are divided by $d$, as in equation (\ref{eq:10}).)

Start with a function $g$ on $Y$ which vanishes to
order $k$ along $E_{v}$.  Viewed on the partial resolution
$Y_1$, $g$ vanishes to the same order in the resolution of the
singular point, which is a splice-quotient for
$\tilde\Gamma$.  Lifting further via $\bar{X_{1}}\rightarrow Y_{1}$,
by induction $g$ is in the $rk$th piece of the filtration
(corresponding to the node $v$ in the splice diagram $\tilde\Delta$) at each of
the $N_{1}$ splice-type singularities; we use $r= \det(\tilde\Gamma)$.
Finally, pulling back leaf variables under the
pseudo-reflection quotient map
$X_{1}\rightarrow \bar{X_{1}}$ multiplies $v$--weights by $d$, via
comparison of $\tilde\Delta$ and $\Delta'$; since the map on the
associated gradeds of corresponding singular points is easily seen to
be injective (simply replace $\bar{A}$ by $A_{1}^{d}$),
the pullback of $g$ to $X_{1}$ now
has weight $drk$ in each of the
$N_{1}$ associated gradeds.

We next look at pulling back a function $h$ from $X$ to $X_{1}$, and
show that the induced map from the $v$--associated graded of $X$ to
the direct sum of the $N_{1}$ associated gradeds on $X_{1}$ multiplies
by degree $r$ and is injective.  Once this claim is established, we
see that the pull-back of our original $g$ to $X$ must have
$v$--filtration weight equal to $dk$, in order to get the correct
weight at each point of $X_{1}$.  The following lemma thus completes
part (\ref{item:main-grading}) of the theorem for these nodes.\qed

\begin{lemma} The map $X_{1}\rightarrow X$ induces a natural map from
    the $v$--associated graded of $X$ to the direct sum
of $N_{1}$ associated gradeds for the splice singularities
on $X_{1}$.  This map is injective and multiplies degrees by $r$.
\end{lemma}
\begin{proof} The map on polynomial rings
  $$R=\mathbb C [x_{i},y_{j}]\rightarrow R_{1}=\mathbb C
  [A_{i},B_{j}]$$
  given by equation (\ref{eq:5}) gives $X_{1}\rightarrow X$,
  which is exactly $rN_{1}$--to-one off the locus $x_{1}=0$. The ideal
  $I\subset R$ generated by the $v$--leading forms of the splice
  equations includes the forms
  $x_{i}^{n_{i}}+a_{i}x_{\sigma}^{n_{\sigma}}$, $1\leq i \leq \sigma
  -1$.  Let $J\subset R_{1}$
  be the ideal generated by the proper transforms of the elements of
  $I$ (ie, factor out by the highest power of $A_{1}$ occurring in
  any equation).  Then the induced map $R/I \rightarrow R_{1}/J$ gives
  a map on spectra which is surjective off the image of the locus
  $x_{1}=0$.  But the $v$--associated graded $R/I$ defines a reduced
  two-dimensional complete intersection, and $x_{1}=0$ is a curve on
  it (Theorem \ref{th:splice is CI}); thus, $R/I \rightarrow R_{1}/J$
  must be injective.

    Next, assign ``$v'$--weights'' to the variables of $R_{1}$:
$$v'(A_{1})=\ell_{vv^{*}}/N; \quad v'(A_{i})=0,~ i>1;\quad
v'(B_{j})=r\ell_{vw_{j}}-\ell'_{v^{*}w_{j}}\ell_{vv^{*}\,},$$
where $w_{j}$ is the leaf corresponding to $y_{j}$ or $B_{j}$.  It is
easy to check that under $R\rightarrow R_{1}$, $v$--weights give
$r$ times $v'$--weights.  The same is true for the graded injection
$R/I \rightarrow R_{1}/J$.

Finally, the variety corresponding to $J$ has the equations
\begin{align*}
1+a_{1}A_{\sigma}^{n_{\sigma}}&=0\\
A_{i}^{n_{i}}+a_{i}A_{\sigma}^{n_{\sigma}}&=0,~ i\geq 2\,,
\end{align*}
which define $N_{1}$ reduced points
$\bar{c}=(c_{2},\ldots,c_{\sigma})$.  Thus, $R_{1}/J$ is a direct sum
of $N_{1}$ graded quotients $R_{1}/J_{\bar{c}}$, where $J_{\bar{c}}$ is
obtained by replacing
$A_{i}$ by $c_{i}$ (for $i>1$) in the defining equations of $J$.  It follows as
in the discussion of \ref{le:tilde} that each such quotient is the
associated graded of the $\Delta'$--splice diagram singularity
corresponding to the node $v$ (or equivalently the weight $v'$) at the
point of $X_{1}$ corresponding to $\bar{c}$, ie, at the point
$A_{1}=B_{j}=0$ (all $j$) and $A_{i}=c_{i}$ ($i\geq 2$).  This
completes the proof of the lemma.
\end{proof}

\section{Naturalness of splice diagram equations}

The definition of splice type equations in Theorem \ref{th:main} might
appear to depend on the choice of monomials satisfying the relevant
conditions, but in fact it does not: for a given node $v$ and edge $e$
at $v$, any two such choices of monomial differ by something of higher
order, which can then be absorbed in the higher order terms of the
splice diagram equation. Precisely:
\begin{theorem}\label{th:change monomial}
  Suppose $M=M_{ve}$ and $M'=M'_{ve}$ are two admissible monomials for
  $\Delta$ that satisfy the $D(\Gamma)$--equivariance condition.  Then
  for some $a\in \C^*$, $M'-aM$ has $v$--weight greater than $d_v$.
  In particular, the corresponding notions of splice diagram
  equations are the same.\comment{minor changes, 7/15}
\end{theorem}
\begin{proof}Choose splice diagram equations as in Theorem
  \ref{th:main}, let $(X,o)$ be the resulting complete
  intersection singularity, and let $Y=X/D(\Gamma)$. Thus $Y$ has a good
  resolution $\bar Y$ with dual graph $\Gamma$. Denote
  $d=\disc(\Gamma)=|D(\Gamma)|$

  An analytic function on $Y$ is simply a $D(\Gamma)$--invariant
  function on $X$, and thus has a $v$--weight for each node of
  $\Delta$.  It also induces a function on the resolution $\bar Y$,
  and thus has an order of vanishing on the exceptional curve $E_v$
  corresponding to $v$.  We recall from item (\ref{item:main-grading})
  of Theorem \ref{th:main} that, for any node
  $v$, the $v$--weight of a function $f$ on $Y$ is $d$ times its order
  of vanishing on the corresponding exceptional curve $E_v$ of the
  resolution.

  Now let $E_1,\dots,
  E_\delta$ be the exceptional curves that intersect $E_v$,
  corresponding to edges $e_1,\dots,e_\delta$ at $v$. Choose an
  admissible monomial $M_i$ that satisfies the
  $D(\Gamma)$--equivariance condition for each edge $e_i$ at $v$. Each
  $M_i^d$ is $D(\Gamma)$--invariant and hence defined on $\bar Y$ and,
  by the above remark, it vanishes to order $d_v$ on the exceptional curve
  $E_v$.  In the same way, at an adjacent node $v^*$ of the splice
  diagram as in
$$\splicediag{8}{30}{
  &&&&\\
  \Vdots&\overtag\Circ v {8pt}\lineto[ul]_(.5){n_1}
  \lineto[dl]^(.5){n_{\delta-1}}
  \lineto[rr]^(.25){n_\delta}^(.75){m_{\nu}}&& \overtag\Circ{v^*}{8pt}
  \lineto[ur]^(.5){m_1}
  \lineto[dr]_(.5){m_{\nu-1}}&\Vdots\\
  &&&&\hbox to 0 pt{~,\hss} }$$
the order of vanishing of $M_i^d$
on $E_{v^*}$ is $d_vd_{v^*}/(n_\delta m_\nu)$ for $i\ne\delta$ and is
$n_\delta m_\nu$ for $i=\delta$.  In particular, on $E_{v^*}$,
$M_\delta^d$ vanishes to order $D$ more than the other $M_i^d$'s,
where $D$ is the edge determinant for the edge $e_\delta$.

In the maximal splice diagram we have a node for every exceptional
curve and all edge determinants are $d$ (Theorem \ref{th:props}). So
we have shown that on each $E_i$ that intersects $E_v$, the
$M_j^d$ with $j\ne i$ vanish to a common order and $M_i^d$ vanishes to
order $d$ greater.  Thus $M_i^d/M_\delta^d$ on $E_v$
has a zero of order $d$ at $E_v\cap E_i$, a pole of order $d$ at
$E_v\cap E_\delta$, and no other zero or pole.

Let the edge $e$ in the theorem be $e=e_1$. Since $M_1$ and $M_\delta$
transform the same way under $D(\Gamma)$, $M_1/M_\delta$ is defined on
$\bar Y$, and $M_1/M_\delta$ on $E_v$ has a simple zero at $E_v\cap
E_i$, a simple pole at $E_v\cap E_\delta$, and no other poles or
zeros. Any other choice $M'_1$ for $M_1$ gives identical zero and pole
for $M'_1/M_\delta$. So for some $a\in \C^*$, $(M'_1 - a
M_1)/M_\delta$ vanishes identically on $E_v$, whence $(M'_1-aM_1)^d$
vanishes to higher order on $E_v$ than does $M_\delta^d$. Since the
$v$--weight of a function $f$ is measured by the order of vanishing of
$f^d$ on $E_v$, the first assertion of the theorem follows.

To prove the second statement of the Theorem, we must show that modulo
the equations defining $(X,o)$, $M'$ is equal to $aM$ plus monomials
of higher $v$--weight. This is the definition of the weight
filtration on $X$ (for convenience, we assume that the defining
equations are polynomials): if $P$ is the graded
polynomial ring in our variables, $I_k$ its ideal generated by
monomials of $v$--weight $k$, and $J$ the defining ideal for $X$, then the
weight filtration on $P/J$ has $k$-th piece $(I_k +J)/J$;  so, modulo
$J$, anything of weight at least $k$ in $P/J$ can be written as a sum
of monomials of degree at least $k$ in the polynomial ring.
\end{proof}

\section{Semigroup and congruence conditions in the two-node case}
\label{sec:twonode}

We first revisit the congruence condition of Proposition
\ref{prop:congruence} in the case that the edge $e$ connects $v$ to
an end node. Thus suppose we are in the following situation, where for
  convenience in this section we will assume minimal good resolutions
  (although this is not essential) and $e$ is the edge of $\Delta$ from $v$ to $v^*$.
$$
\xymatrix@R=3pt@C=24pt@M=0pt@W=0pt@H=0pt{\\
&\lefttag{\Circ}{n_\sigma/p_\sigma}{8pt}\dashto[ddrr]&
&&&&\\&&&&&&&&&& \\&
&&\Circ\lineto[dr]&&{\hbox to 0 pt {$\scriptstyle n/p~\rightarrow$\hss}}
&&&&\Circ\lineto[dl]\dashto[ur]\\
\Gamma=\quad&
&\hbox to 0pt{\hglue-8pt\vbox to 0 pt{\vss\vss\vdots\vss}\hss}&&
\undertag{\overtag{\Circ}{-b}{8pt}}{v^*}{6pt}\lineto[r]&
\dashto[rr]&&\lineto[r]&\undertag\Circ{v}{6pt}&&\hbox to 0
pt{\hss\vbox to 0 pt{\vss\vss\vdots\vss}\hglue5pt}
\\&&&\Circ\lineto[ur]&&&&&&\Circ\lineto[ul]\dashto[dr]\\
&&&&&&&&&&
\\&\lefttag{\Circ}{n_1/p_1}{8pt}\dashto[uurr]
}
$$$$
\xymatrix@R=6pt@C=30pt@M=0pt@W=0pt@H=0pt{\\
&\Circ\lineto[ddrr]^(.7){n_\sigma}&&&&&&&\\&&&&&&&\dashto[ur]
\\
\Delta=&&\hbox to 0pt{\hglue-8pt\vbox to 0 pt{\vss\vss\vdots\vss}\hss}
&\undertag\Circ{v^*}{6pt}\lineto[rrr]^(.2)r^(.8)s&&&
\undertag\Circ{v}{6pt}\lineto[ur]^(.6){m_\mu}\lineto[dr]_(.6){m_1}
&\hbox to 0 pt{\hss\vbox to 0 pt{\vss\vss\vdots\vss}\hglue-8pt}&
\\&&&&&&&\dashto[dr]\\&\Circ\lineto[uurr]_(.7){n_1}&&&&&&&
\\&}$$
 As usual, we
represent strings in $\Gamma$ by their continued fractions (the
continued fraction for the empty
string is $n/p=1/0$).
Denote $N=\prod_{i=1}^\sigma n_i$, $M=\prod_{j=1}^\mu m_j$, $d=\disc(\Gamma)$.
Then,
using Proposition \ref{prop:endwt}, the condition of Proposition
\ref{prop:congruence} for an admissible monomial $M_{ve}=\prod
x_i^{\alpha_i}$ is that for $i=1,\dots,\sigma$
$$\Bigl[\sum_{j\ne
  i}\alpha_jNr/(n_in_jd)+\alpha_i(Nr/(n_i^2d)+p'_i/n_i)\Bigr]=[MN/(n_id)]\,.$$
This simplifies to
$$[rs/(n_id)+\alpha_ip'_i/n_i]=[MN/(n_id)]\,,$$ 
$$[(rs-MN)/(n_id)+\alpha_ip'_i/n_i]=[0]\,.\leqno{\rm or}$$
Since, by Proposition \ref{prop:edge det},
$rs-MN=dn$,
the above is equivalent to
$$p'_{i}\alpha_{i} \equiv -n\quad
(\text{mod }
n_{i})\,,\leqno{\rm or equivalently}$$
$$\alpha_{i}\equiv -np_{i}\quad (\text{mod }n_{i})\,.$$
Now, solutions
of these congruences may be written
\begin{align*}
  \alpha_{i}&=n_{i}\lceil np_{i}/n_{i}\rceil-np_{i}+n_{i}\delta_{i}\,,
\end{align*}
where $\lceil x\rceil$
means least integer $\geq x$; further the non-negativity of the
$\alpha_{i}$ is equivalent to the non-negativity of
the $\delta_{i}$. Thus the equality $s=\sum\alpha_iN/n_i$, which
expresses that the monomial $M_{ve}$ is admissible, can be written
\begin{equation}
\label{eq:2.1}
s=\sum_i\left(N\lceil np_i/n_i\rceil-Nnp_i/n_i+N\delta_i\right)\,.
\end{equation}
On the other hand, by computing determinant of a
star-shaped graph  we get
\begin{equation*}
s=Nn(b-\sum_{i=1}^{\sigma}p_{i}/n_{i}\ -p/n)\,.
\end{equation*}
Thus
formula (\ref{eq:2.1})  is
equivalent to
\begin{align*}
\sum \delta_{i} &= nb-p-\sum \lceil np_{i}/n_{i}\rceil
\end{align*}
If the right hand side of this expression is non-negative, then
non-negative $\delta_{i}$ can be found so that the corresponding
$\alpha_{i}$ satisfy both the congruence conditions \emph{and} the
semigroup condition at the given edge.  We summarize in the
\begin{proposition}\label{prop:was two node}
  Consider the edge $e$ leading from $v$ to an end-node as above.
  Then the following inequality is necessary and sufficient in order
  that the semigroup and congruence conditions are both satisfied for
  vertex $v$ and edge $e$:
  $$nb-p-\sum\lceil np_{i}/n_{i}\rceil \geq 0\,.\eqno{\qed}$$
\end{proposition}
This proposition has the immediate corollary:
\begin{proposition}
  \label{prop:two node}
  The following two-node resolution graph (dashed lines represent strings
  described by continued fractions starting from the interior\break weights;
  the central string $n/p$ starts from the left node)
$$
\xymatrix@R=6pt@C=24pt@M=0pt@W=0pt@H=0pt{
\Circ&& &&& &&\Circ\\
\hbox to 0 pt{\hss$\Gamma=$\qquad}&\Vdots&\overtag\Circ{-b}{9pt}
\dashto[ull]_(.5){n_1/p_1}
\dashto[dll]^(.5){n_\sigma/p_\sigma}
\dashto[rrr]^(.45){n/p~\rightarrow}&&&
\overtag\Circ{-c}{9pt}
\dashto[urr]^(.5){m_1/q_1}
\dashto[drr]_(.5){m_\tau/q_\tau}&\Vdots\\
\Circ&& &&& &&\Circ
}
$$
satisfies the semigroup and congruence conditions if and only if
$$
  nb-p-\sum\lceil np_{i}/n_{i}\rceil \geq 0 $$$$
  nc-p'-\sum\lceil nq_{j}/m_{j}\rceil \geq 0 \,.\eqno{\qed}
$$
\end{proposition}
We remark that the negative-definiteness of this graph $\Gamma$ is
equivalent to the condition that the edge determinant
$rs-\prod_in_i\prod_jm_j$ is positive together with the positivity of
$s$ and $r$. The latter are slightly weaker conditions than those
of the proposition; $s/N>0$ and $r/M>0$ can be written:
\begin{gather*}
s/N=nb-p-\sum np_{i}/n_{i} > 0\\
r/M=nc-p'-\sum nq_{j}/m_{j} > 0\,.
\end{gather*}
According to our main theorem, $\Gamma$ occurs as the resolution dual
graph of a splice-quotient singularity, that is, a surface singularity
whose universal abelian cover is of splice type, if the semigroup and
congruence conditions are satisfied.  We had conjectured earlier that
a $\Q$--Gorenstein singularity with $\Q$HS link is
always of this type. Although counter-examples are now known (see
\cite{nemethi et al}), the conjecture appears to hold in a surprising
number of cases. Singularities with rational or minimally elliptic
resolution graphs are automatically $\Q-$Gorenstein (even Gorenstein
for minimally elliptic, \cite{laufermin}), and we assert the truth of
the following important

\begin{conjecture}\label{rational}  Let $Y$ be a rational or
$\Q$HS-link minimally
  elliptic surface singularity.  Then $Y$ is a splice-quotient
  singularity as in \ref{splicequotient}.
\end{conjecture}

\noindent{\bf Note added April 2005}\qua T Okuma \cite{okuma} has
announced a proof of this conjecture, see Section \ref{sec:okuma}. In
the original version of this paper we gave here a partial proof of
this conjecture in the two-node case, which we now omit.

\section{Appendix 1: Splicing and plumbing}\label{sec:splicing}

This appendix reviews in more detail how a splice diagram is
associated to a resolution diagram and explains why it is a
topological invariant of the $3$--manifold link. We also prove some
technical results needed earlier in the paper.

Recall (see the beginning of Section \ref{sec:basics}) that a
\emph{splice diagram} is a finite tree with no valence $2$ vertices,
decorated with integer weights as follows: for each node $v$ and edge
$e$ incident at $v$ an integer weight $d_{ve}$ is given. Thus an edge
joining two nodes has weights associated to each end, while an edge
from a node to a leaf has just one weight at the node end.  Moreover,
we will show that the splice diagrams which arise in the study of
links of singularities always satisfy the following conditions:
\begin{itemize}
\item All weights are positive.
\item All edge determinants are positive.
\item The ideal condition (Definition \ref{def:ideal}).
\end{itemize}
(For the splice diagram
associated with an arbitrary graph-manifold rational homology sphere
the first two conditions need not hold.) In the process, we will also
have need for a variant of splice diagrams where valency 2 vertices
are permitted, and weights are also associated to the leaf end of an
edge ending in a leaf.

In \cite{eisenbud-neumann} splice diagrams were used (among other
things) to classify the topology of integral homology sphere
singularity links.  The splice diagrams that arise this way are
precisely the splice diagrams as above with pairwise coprime positive
weights around each node (in which case the ideal condition is
automatic).  \comment{is this a silly insertion?}  The paper
\cite{neumann-wahl02} was the first to associate a splice diagram more
generally to any \emph{rational} homology sphere singularity link
$\Sigma$. The splice diagram no longer determines the topology of
$\Sigma$, but we claim that it does determine the topology of the
universal abelian cover of $\Sigma$ (which is $\Sigma$ itself, if
$\Sigma$ is a $\Z$--homology sphere).  The current paper establishes
this assertion only when the semigroup and congruence conditions
(\ref{def:sg} and \ref{def:congruence}) are satisfied, but the result
holds without these conditions, and extends even to arbitrary
graph-manifold homology spheres; this will be proved elsewhere.

For ease of exposition we restrict to the singularity link case here.
In this case we can describe the splice diagram in terms of a
resolution of the singularity. This was described briefly at the start
of Section \ref{sec:congruence condition} but we will recall it in
more detail.

Thus, let $(Y,o)$ be a normal surface singularity germ and $\Sigma$
its link, that is, the boundary of a regular neighborhood of $o$ in
$Y$.  Assume that $\Sigma$ is a rational homology sphere,
equivalently, $H_1(\Sigma)$ is finite.  Let $\pi\colon\bar Y\to
Y$ be a good resolution. ``Good'' means that the exceptional divisor
$E=\pi^{-1}(o)$ has only normal crossings.  The rational homology
sphere condition is equivalent to the condition that $E$ is
\emph{rationally contractible}; that is,
\begin{itemize}
\item each component of $E$ is a smooth rational curve;
\item the dual resolution graph $\resgraph$ (the graph with a vertex for
  each component of $E$ and an edge for each intersection of two
  components) is a tree.
\end{itemize}

We weight each vertex $v$ of $\resgraph$ by the self-intersection
number $E_v\cdot E_v$ of the corresponding component $E_v$ of $E$. The
\emph{intersection matrix} for $\resgraph$ is the matrix $A(\resgraph)$
with entries $a_{vw}=E_v\cdot E_w$.
It is well known that $A(\resgraph)$ is negative-definite and its
cokernel (also called the \emph{discriminant
group}) is $H_1(\Sigma)$.
In particular, $\disc(\resgraph):=\det(-A(\resgraph))$
is the order of $H_1(\Sigma)$.

A \emph{string} in $\resgraph$ is a connected subgraph consisting of
vertices that have valency $\le2$ in $\resgraph$. The resolution is
\emph{minimal} if no $(-1)$--weighted vertex of $\Gamma$ occurs on a
string. We do not necessarily want to assume minimality here.

The splice diagram $\Delta$ for $\Sigma$ has the same overall shape as
the resolution graph $\resgraph$; it's underlying graph is obtained from
$\resgraph$ by suppressing valency two vertices.  The weights on edges are
computed by the following procedure: At a vertex $v$ of $\Delta$ let
$\resgraph_{ve}$ be the subgraph of $\resgraph$ cut off by the edge of
$\resgraph$ at $v$ in the direction of $e$, as in the following picture.
The corresponding weight is then $d_{ve}:=\disc(\resgraph_{ve})$.
$$
\xymatrix@R=6pt@C=24pt@M=0pt@W=0pt@H=0pt{
&&&&&&\\
\Vdots&&\undertag{\overtag{\Circ}{a_{vv}}{8pt}}{v}{8pt}
\lineto[rr]\dashto[ull]\dashto[dll]&{}\undertag{}{e}{8pt}&
\overtag{\Circ}{a_{ww}}{8pt}\dashto[urr]\dashto[drr]&&\Vdots\\
&&&&&&\\
&&&&&{\hbox to 0pt{\hss$\underbrace{\hbox to 70pt{}}$\hss}}&\\
&&&&&{\resgraph_{ve}}&}$$

\begin{example}\label{ex:2}
  Here is an example of a resolution graph with integral homology
  sphere link.  The reader can check that $A(\resgraph)$ is
  negative-definite and unimodular (a quick method is given in
  \cite{eisenbud-neumann}).
$$
\xymatrix@R=6pt@C=24pt@M=0pt@W=0pt@H=0pt{
\\
&\overtag{\Circ}{-2}{8pt}&&&&\overtag{\Circ}{-2}{8pt}\\
{\resgraph\quad=}&&\overtag{\Circ}{-1}{8pt}\lineto[ul]\lineto[dl]\lineto[r]&
\overtag{\Circ}{-17}{8pt}&\overtag{\Circ}{-1}{8pt}\lineto[ur]\lineto[dr]\lineto[l]&\\
&\overtag{\Circ}{-3}{8pt}&&&&\overtag{\Circ}{-3}{8pt}\lineto[r]&\overtag{\Circ}{-2}{8pt}}
$$
Its splice diagram is:
$$
\splicediag{6}{30}{
&\Circ&&&\Circ\\
\Delta\quad=&&\Circ\lineto[ul]_(.25){2}\lineto[dl]^(.25)3
&\Circ\lineto[dr]_(.25){5}\lineto[ur]^(.25)2
\lineto[l]_(.2){11}_(.8){7}\\
&\Circ&&&\Circ
}\qquad$$
For example, the weight $7$ on the left node of $\Delta$
is $\disc(\resgraph_{ve})$ with
$$
\xymatrix@R=6pt@C=24pt@M=0pt@W=0pt@H=0pt{
&&&\overtag{\Circ}{-2}{8pt}\\
{\resgraph_{ve}}\quad=&\overtag{\Circ}{-17}{8pt}\lineto[r]&
\overtag{\Circ}{-1}{8pt}\lineto[ur]\lineto[dr]\\
&&&\overtag{\Circ}{-3}{8pt}\lineto[r]&\overtag{\Circ}{-2}{8pt}}
$$
Here is another resolution graph with the same splice diagram
$$
\xymatrix@R=6pt@C=24pt@M=0pt@W=0pt@H=0pt{
\\
&\overtag{\Circ}{-2}{8pt}&&&\overtag{\Circ}{-2}{8pt}\\
&&\overtag{\Circ}{-3}{8pt}\lineto[ul]\lineto[dl]\lineto[r]
&\overtag{\Circ}{-2}{8pt}\lineto[ur]\lineto[dr]\lineto[l]&\\
\overtag\Circ{-2}{8pt}\lineto[r]&\overtag{\Circ}{-2}{8pt}&&&
\overtag{\Circ}{-2}{8pt}\lineto[r]&\overtag{\Circ}{-2}{8pt}\lineto[r]
&\overtag{\Circ}{-2}{8pt}\lineto[r]&\overtag{\Circ}{-2}{8pt}\hbox to 0
pt{~.\hss}\\&}
$$
It has discriminant 17, so its link has first homology $\Z/17$.
\end{example}

If $\Sigma$ is a $\Z$--homology sphere, then the minimal resolution
graph can be recovered from the splice diagram; an algorithm to do
this is is described in \cite{neumann-wahl-zsphere}, improving on a
procedure in \cite{eisenbud-neumann}. Thus, in the above example,
$\Gamma$ is the only minimal resolution graph with splice diagram
$\Delta$ and with $\Z$--homology sphere link. However, there can be
several minimal resolution graphs with the same splice diagram
representing $\Q$--homology spheres (infinitely many if the splice
diagram has just one node and finitely many otherwise).

To understand the resolution graphs that correspond to a given
splice diagram it is helpful to consider the
\emph{maximal splice diagram}:
the version of the splice diagram that we get from the resolution
graph if we do not first
eliminate vertices of valency $2$, and include edge weights at
\emph{all} vertices --- also the leaves.  Thus, for the first of
Examples \ref{ex:2}, the maximal splice diagram is:
$$
\xymatrix@R=6pt@C=24pt@M=0pt@W=0pt@H=0pt{
  &\Circ&&&&\Circ\\
  {\Delta'\quad=}&&\Circ\lineto[ul]_(.25){2}_(.75){11}
\lineto[dl]^(.25){3}^(.75)5
\lineto[r]^(.25)7^(.75)1&
\Circ&\Circ\lineto[ur]^(.25)2^(.75){28}\lineto[dr]_(.25)5_(.75)9
\lineto[l]_(.25){11}_(.75)1&\\
&\Circ&&&&\Circ\lineto[r]_(.25)2_(.75)5&\Circ}
$$
and for the second it is
$$
\xymatrix@R=6pt@C=24pt@M=0pt@W=0pt@H=0pt{
&\Circ&&&\Circ\\
&&\Circ\lineto[ul]_(.25)2_(.75){19}\lineto[dl]^(.25)3^(.75){15}
\lineto[r]^(.25)7^(.75){11}
&\Circ\lineto[ur]^(.25){2}^(.75){36}\lineto[dr]_(.25)5_(.75){21}
\\
\Circ\lineto[r]^(.25){16}^(.75)2&\Circ&&&
\Circ\lineto[r]^(.25)4^(.75){20}&\Circ\lineto[r]^(.25)3^(.75){19}
&\Circ\lineto[r]^(.25)2^(.75){18}&\Circ\hbox to 0
pt{~.\hss}\\&}
$$
The maximal splice diagram has the following properties.
\begin{theorem}\label{th:props}
  {\rm(1)}\qua For any pair of vertices $v$ and $w$ of the maximal diagram
  let $\ell_{vw}$ be the product of the weights adjacent to, but not on,
  the shortest path from $v$ to $w$ in $\Delta'$ (in particular,
  $\ell_{vv}=d_v$, the product of weights at $v$). Then the matrix
  $L:=(\ell_{vw})$ satisfies $\frac1{\disc(\resgraph)}L=-A(\resgraph)^{-1}$.

  {\rm(2)}\qua Every edge determinant for the maximal splice diagram
  is $\disc(\resgraph)$.
 \end{theorem}
\begin{proof}
  Property (1) of the theorem says \begin{equation}\label{eq:1}
A(\resgraph)L=-\disc(\resgraph)I\,,
\end{equation}
which is easily shown by
computing that the adjoint matrix of $-A(\resgraph)$
equals $L$.  This calculation is carried out explicitly in
Lemma 20.2 of \cite{eisenbud-neumann}.

For property (2), suppose we have an edge connecting vertices $v$ and
$w$ of the maximal splice diagram as follows, $$\splicediag{8}{40}{
  &&&\\&&&  \\
  \Vdots&\undertag\Circ
  v{4pt}\lineto[uul]_(.35){n_1}\lineto[ul]^(.45){n_2}\lineto[dl]^(.35)
  {n_\sigma} \lineto[r]^(.25){r}^(.75){s}&\undertag\Circ w{4pt}
  \lineto[uur]^(.35){m_1}\lineto[ur]_(.45){m_2}\lineto[dr]_(.35){m_\tau}
  &\Vdots\\
  &&&\\& }$$
and write $N=\prod_1^\sigma n_i$, $M=\prod_1^\tau m_j$,
$N_i=N/n_i$ (if $v$ or $w$ is a leaf the corresponding $N$ or $M$ is
$1$).  For each $i=1,\dots,\sigma$ let $L_i$ be the product of the
weights just beyond the other end of the $n_i$--weighted edge.  Then
the $vv$-- and $vw$--entries of equation (\ref{eq:1}) are:
\begin{align*}
Nra_{vv}+NM+\sum_1^\sigma rN_iL_i&=-\disc(\resgraph)\\
NMa_{vv}+sM+\sum_1^\sigma MN_iL_i&=0\,.
\end{align*}
Multiplying the second of these equations by $r/M$ and then
subtracting the first from it gives the desired equation
$rs-MN=\disc(\resgraph)$.
\end{proof}

We can generalize part 2 of the above theorem to any edge of a splice
diagram.
Let $\resgraph$ be a resolution graph and $\Delta$ its splice diagram. Thus
each edge of $\Delta$ corresponds to a string in $\resgraph$.
\begin{proposition}\label{prop:edge det}
  Let $e$ be an edge of $\Delta$ corresponding to a string $E$ of
  $\resgraph$. Then the edge determinant $D(e)$ is given by
  $$D(e)=\disc(E)\disc(\resgraph)\,,$$
where $\disc(E)=1$ if $E$ is the empty string.
\end{proposition}
\begin{proof}
We need some preparation.
\begin{lemma}\label{le:previous}
  Suppose in $\Gamma$ we have an extremal string with
  continued fraction $n/p=b_1-1/\dots-1/b_k$, and associated splice
  diagram as follows:
\begin{gather*}
\xymatrix@R=4pt@C=24pt@M=0pt@W=0pt@H=0pt{
&\dashto[ddrr]\\ \\
\Gamma=\quad
&\hbox to 0pt{\hglue0pt\vbox to 0 pt{\vss\vss\vdots\vss}\hss}&&
\undertag{\overtag{\Circ}{-b}{8pt}}{v}{4pt}\lineto[r]&
\overtag\Circ{-b_1}{8pt}
\dashto[rr]&&\overtag\Circ{-b_{k-1}}{8pt}
\lineto[r]&\overtag\Circ{-b_k}{8pt}\\
&&&&&\vbox to 0 pt{\vss\hbox to 0 pt{\hss\hss$\scriptstyle
n/p~\rightarrow$\hss}\vss}
\\&\dashto[uurr]
}
\\
\xymatrix@R=6pt@C=20pt@M=0pt@W=0pt@H=0pt{\\
&\lineto[ddrr]^(.6){n_\sigma}&&&&&&&\\
\\
\quad\Delta=&&\hbox to 0pt{\hglue-8pt\vbox to 0 pt{\vss\vss\vdots\vss}\hss}
&\undertag\Circ{v}{4pt}\lineto[rrr]^(.2)n&&&\Circ
\\ \\&\lineto[uurr]_(.6){n_1}&&&&&&&
\\&}
\end{gather*}
Let $\Gamma_0$ be the result
of removing the string, so it
consists of $v$ and what is to the left.
Then, with $N=n_1\dots n_\sigma$,  $$\disc(\Gamma)=n\disc(\Gamma_0)-Np\,.$$
\end{lemma}
\begin{proof}
  This is the edge determinant equation of part (2) of Theorem
  \ref{th:props} applied to the edge from $v$ to the $-b_1$--weighted
  vertex in the maximal splice diagram, since the determinant of the
  string starting at $-b_2$ is $p$.
\end{proof}
The following lemma has been used earlier (eg, Proposition
\ref{prop:endwt}), since, even though the edge weight at a leaf is not
part of the data of a splice diagram, it is needed in computing
discriminant groups.
\begin{lemma}\label{le:end weight}
Suppose we have a leaf $w$ of a
splice diagram
$$
\xymatrix@R=6pt@C=36pt@M=0pt@W=0pt@H=0pt{&&\\
\Delta=&\Vdots&\undertag{\Circ}{v}{4pt}\lineto[ul]_(.4){n_1}
\lineto[dl]^(.4){n_\rho}
\lineto[r]^(.3){n}&\undertag{\Circ}{w}{4pt}\\&& \\&}
$$
resulting from a resolution string with continued fraction
$n/p=b_1-1/\dots-1/b_k$.  Denote by $p'$ the discriminant of
the\comment{06/25 regression; original version was better}
string with $w$ removed, so $pp'\equiv 1$ (Lemma \ref{le:cf}; if
the string is \minimal{} then $p'$ is the unique such positive integer
with $p'\le n$). Denote $N=n_1\dots n_\sigma$.
  Then the splice diagram weight $x$ at $w$ is given by
  $$x=\frac{p'}{n}\disc(\resgraph)+\frac{N}{n}$$
\end{lemma}
\begin{proof}
  Denote $n'=(pp'-1)/n$, so $pp'-nn'=1$.  By Lemma
  \ref{le:cf} we have $p'/n'=b_1-1/\cdots-1/b_{k-1}$ (if $k=1$, then
  $p=p'=1$, $n'=0$).  Apply the previous Lemma \ref{le:previous} to
  both $\Gamma$ and the result $\Gamma'$ of removing the rightmost
  vertex of $\Gamma$.  This gives the following (also if $k=1$):
$$\disc(\Gamma)=n\disc(\Gamma_0)-Np,\qquad
x:=\disc(\Gamma')=p'\disc(\Gamma_0)-Nn'\,.$$
Solving the first of these equations for
$\disc(\Gamma_0)$ and inserting in the second gives
$$x=p'\left(\frac{\disc(\Gamma)}{n}+\frac{Np}n\right)-Nn'=
\frac{p'}n\disc(\Gamma)+\frac Nn\,,$$
as desired.
\end{proof}

We now complete the proof of Proposition \ref{prop:edge det} by
induction on the length of the string. We already know it for the
empty string by Theorem \ref{th:props}, so suppose we have partially
reduced the maximal splice diagram as follows:
$$
\xymatrix@R=6pt@C=24pt@M=0pt@W=0pt@H=0pt{&&&&&&&\\
\Vdots&&\undertag\Circ{v}{4pt}
\lineto[ull]_(.3){n_1}\lineto[dll]^(.3){n_\sigma}\lineto[rr]^(.2){r}^(.8){s_1}
&&\undertag\Circ{v_1}{4pt}\lineto[r]^(.3){r_1}^(.7){s_2}&
\undertag\Circ{v_2}{4pt}
\lineto[urr]^(.3){p_1}\lineto[drr]_(.3){p_\rho}&&\Vdots&(\rho\ge0)\,,\\
&&&&&&&\\~}
$$
where $v_1$ and $v_2$ were adjacent in the maximal splice diagram, but
the edge from $v$ to $v_1$ may correspond to a non-empty string. Denote the
string of $\Gamma$ between $v$ and $v_i$ by $E_i$, $i=1,2$, and
denote $N=n_1\dots,n_\sigma$, $P=p_1\dots,p_\rho$.

By Theorem \ref{th:props} and by the induction assumption we have
$$r_1s_2-s_1P=\disc(\Gamma),\qquad
rs_1-Nr_1=\disc(E_1)\disc(\Gamma)\,.$$
Multiplying the first of these by $N$ and the second by $s_2$ and
adding gives
$$s_1(rs_2-NP)=\disc(\Gamma)(s_2\disc(E_1)+N)\,.$$
Apply the last
lemma (Lemma \ref{le:end weight}) to the result of deleting from
$\Gamma$ the vertex $v_2$ and all
to the right of it. This has determinant $s_2$, so  Lemma \ref{le:end
  weight} gives
$$s_1=\frac{\disc(E_1)}{\disc(E_2)}s_2+\frac N{\disc(E_2)}\,.$$
Inserting this in the previous equation and simplifying gives
$$rs_2-NP=\disc(\Gamma)\disc(E_2)\,,$$
completing the inductive step.
\end{proof}

We earlier needed to understand what happens to weights in a
splice diagram $\Delta$ when part of the resolution diagram $\resgraph$
changes. Of course weights only change if given by determinants of
parts of the resolution graph that have changed. A typical situation
might be the following:
$$
\xymatrix@R=6pt@C=24pt@M=0pt@W=0pt@H=0pt{&&\\
  \Delta=& &\undertag{\Circ}{v}{5pt}\dashto[ul]\dashto[dl]
  \lineto[rr]^(.25){a}^(.75){b}&&
  \undertag{\Circ}{w}{5pt}\dashto[ddl]\dashto[ddr]
  \lineto[rr]^(.25){a'}&&\dashto[r]&\frame[5pt]{\Delta_0}\\&& \\
  &&&&&\\&}$$
where changing the part of the resolution diagram
$\resgraph$ corresponding to $\Delta_0$ will not change $b$ but will
change the weights $a$ and $a'$, say to $\tilde a$ and $\tilde
a'$. Denote the changed resolution diagram by $\tilde \resgraph$ and
the corresponding splice diagram by $\tilde \Delta$
\begin{lemma}
  Let $M$ be the product of weights other than $a$ at $v$ and $L$
  the product of weights other than $b$ and $a'$ at $w$.
$$a\disc(\tilde\resgraph)-\tilde
a\disc(\resgraph)=\frac{ML}b(a'\disc(\tilde
\resgraph)-\tilde a'\disc(\resgraph))$$
\end{lemma}
\begin{proof}
  Applying Proposition \ref{prop:edge det} to the edge from $v$ to $w$
  in $\Delta$ and $\tilde \Delta$ gives the equations
\begin{align*}
  ab-a'ML&=n\disc(\resgraph)\\
  \tilde ab-\tilde a'ML&=n\disc(\tilde\resgraph)
\end{align*}
where $n$ is the determinant of the resolution string for the given
edge. Eliminating $n$ from these two equations gives the desired equation.
\end{proof}
\begin{lemma}\label{le:10.7}
Suppose that $\resgraph$ is a resolution diagram with a string
$$
\xymatrix@R=6pt@C=24pt@M=0pt@W=0pt@H=0pt{&&&&&&&&\\
\resgraph=&&\overtag{\Circ}{-b_0}{8pt}\dashto[ul]\dashto[dl]\lineto[r]&
\overtag{\Circ}{-b_1}{8pt}\dashto[r]&&\dashto[r]&
\overtag{\Circ}{-b_n}{8pt}\lineto[r]&
\undertag{\overtag{\Circ}{-b_{n+1}}{8pt}}{v^*}{4pt}\dashto[ur]\dashto[dr]\\
&&&&&&&&\\&}
$$
and $\tilde\resgraph$ results by deleting the node $v^*$ at one end of the
string and all beyond it:
$$
\xymatrix@R=6pt@C=24pt@M=0pt@W=0pt@H=0pt{&&&&&&&&\\
  \tilde\resgraph=&&\overtag{\Circ}{-b_0}{8pt}\dashto[ul]\dashto[dl]\lineto[r]&
  \overtag{\Circ}{-b_1}{8pt}\dashto[r]&&\dashto[r]&\overtag{\Circ}{-b_n}{8pt}
&.\\
  &&&&&&&&\\&}
$$
Suppose $\tilde\Delta$ is the corresponding splice diagram, and $a$
and $\tilde a$ are the $\Delta$-- and $\tilde \Delta$--weights towards
$v^*$ at a node $v$ of\/ $\tilde \Gamma$. Then
  $$a\disc(\tilde\resgraph)-\tilde
  a\disc(\resgraph)=MN(\ell'_{v^*v})^2\,,$$
where
  $M$ is the product of weights of $\Delta$ other than $a$ at $v$
  and $N$ is the product of weights of $\Delta$ at $v^*$ other than the
  weight towards $v$.\end{lemma}
\begin{proof}
  The result is by induction over the distance from $v^*$ to $v$ in
  $\Delta$. The induction step is the previous lemma. The induction
  start is the case that $v^*$ and $v$ are adjacent in $\Delta$. In
  this case the equation to be proved can be written
  $a\disc(\tilde\Gamma)-MN=\tilde a\disc(\Gamma)$. This is the edge
  determinant formula of Proposition \ref{prop:edge det}, since
  $\tilde a$ is the determinant of the string connecting $v$ to $v^*$
  in $\Gamma$ and $\disc(\tilde\Gamma)$ is the edge weight of $\Delta$
  at $v^*$ towards $v$.
\end{proof}

\subsection{Topological description of the splice diagram and Ideal
  Condition}

The weights in a splice diagram have a simple topological meaning.
The standard plumbing description (see, eg, \cite{neumann81}) of the
manifold $\Sigma=\Sigma(\Gamma)$ associated to a resolution graph (or
more general rational plumbing graph) $\Gamma$ shows that to each
string in the graph $\Gamma$ is associated an embedded torus in
$\Sigma$ such that, if one cuts along these tori, $\Sigma$ decomposes
into pieces associated to the leaves and nodes of $\Gamma$. The piece
for a leaf is a solid torus, and for a node is of the form (punctured
disc)$\times S^1$. (If one omits the tori corresponding to leaves,
this essentially describes the JSJ decomposition of $\Sigma$.)  In
particular, the pieces at nodes have natural circle fibers,
topologically determined up to isotopy.

Suppose we have a resolution
diagram and associated splice diagram as follows:
$$\splicediag{8}{30}{
  &&&&&&\frame{\resgraph_1}\\
  \resgraph=& \Vdots&\overtag\Circ{-b_0}{8pt}\dashto[ul]\dashto[dl]
  \lineto[r]&\overtag\Circ{-b_1}{8pt}\dashto[r]&
\dashto[r]&\overtag\Circ{-b_{n+1}~}{9pt}
  \lineto[ur]\lineto[dr]
  &\Vdots\\
  &&&&&&\frame{\resgraph_n} \\
&&&&&\hbox to 0pt{\hss$\underbrace{\hbox to 126pt{}}$\hbox to 12pt{}\hss}\\
&&&&&\hbox to 12pt{\hss$\resgraph_0$\hbox to 9pt{}\hss}\\}$$
$$\splicediag{8}{40}{
  & &&&\\  &&&&\\
  \Delta=& \Vdots&\undertag\Circ v{4pt}\dashto[ul]\dashto[dl]
  \lineto[r]^(.25){d_{ve}}^(.75){d_{v'e}}&\undertag\Circ {v'}{4pt}
  \lineto[uur]^(.35){d_1}\lineto[ur]_(.45){d_2}\lineto[dr]_(.35){d_n}
  &\Vdots\\
  & &&&\\& }$$
(we are denoting $d_i:=d_{v'e_i}$ for the $i$-th edge $e_i$ departing $v'$
to the right).

The topological interpretation of $d_{ve}=\disc(\Gamma_0)$ is simply
that it is the size of $H_1(\Sigma(\Gamma_0);\Z)$; topologically
$\Sigma(\Gamma_0)$ is the manifold one obtains from the right hand
piece after cutting $\Sigma$ along the torus for the edge $e$, by
gluing a solid torus into the boundary torus to kill fibers associated
to the left node $v$ (ie, match them with meridians of the solid
torus).

We can use this for a topological proof that a splice diagram
satisfies the ideal condition of Definition \ref{def:ideal}. Recall
that it says that for each node $v$ and adjacent edge $e$ of a splice
diagram $\Delta$, the edge-weight $d_{ve}$ is in the ideal
$$( \ell'_{vw}:w \text{ a leaf of $\Delta$ in }
  \Delta_{ve} )~\subset~ \Z\,. $$

\begin{definition}
  We will call the positive generator $\bar d_{ve}$ of the above ideal
  the \emph{ideal generator} for $e$ at $v$. So the ideal condition says
  $\bar d_{ve}$ divides $d_{ve}$.
\end{definition}

We refer again to the above diagrams and note that
 $d_i =\lvert
H_1(\Sigma(\resgraph_i);\Z)\rvert$ for $i=0,\dots,n$.  Each of the
manifolds $\Sigma(\resgraph_i)$ contains a knot $K_i$ corresponding to
the edge that attaches $\resgraph_i$ to the rest of $\resgraph$. Note
that the map $H_1(\Sigma(\resgraph_0);\Z)\to
H_1(\Sigma(\resgraph_0)/K_0;\Z)$ is surjective, so $\lvert
H_1(\Sigma(\resgraph_0)/K_0;\Z)\rvert$ divides $d_0$. The
following theorem thus implies the ideal condition.
\begin{theorem}\label{th:topological divisor weight}
  The ideal generator  $\bar d_{ve}$ is $\lvert
  H_1(\Sigma(\resgraph_0)/K_0;\Z)\rvert$.
\end{theorem}
We will prove this theorem inductively, so we
first describe an inductive computation of the ideal generators.
\begin{lemma}\label{le:divisor weights}
  If $v'$ is a leaf ($n=0$) put $\bar d_{ve}=1$.  Inductively, if the
  ideal generator $\bar d_{i}$ is known at the $i$-th edge departing
  $v'$ to the right for each $i$ then $\bar d_{ve}$ is computed as
  $$\bar d_{ve}=\gcd\nolimits_{i=1}^n \Biggl(\bar d_i\,\prod_{j\ne i}^n
d_j\Biggr)\,.$$
\end{lemma}
\begin{proof}
  $\bar d_i\,\prod_{j\ne i}^n d_j$ is the generator of the ideal $(
  \ell'_{v'w}:w \text{ a leaf of $\Delta$ in } \Delta_{v'e_i} ) $, so
  $\gcd_{i=1}^n \bar d_i\,\prod_{j\ne i}^n d_j$ is the generator of
  the ideal $( \ell'_{vw}:w \text{ a leaf of $\Delta$ in } \Delta_{ve}
  )$.
\end{proof}

\begin{proof}[Proof of Theorem \ref{th:topological divisor weight}]
  Define for the moment $\bar d_i=\lvert
  H_1(\Sigma(\resgraph_i)/K_i;\Z)\rvert$ for each $i$. We will show
  that these numbers satisfy the inductive formula of the lemma, so
  they are the ideal generators.

 Let $\Delta_i$ be the subgraph of the splice diagram $\Delta$
  corresponding to the subgraph $\resgraph_i$ of $\resgraph$.
  $\Sigma=\Sigma(\resgraph)$ contains tori $T_i$ corresponding to the
  edges of $\Delta$ that cut off the subdiagrams $\Delta_i$. The torus
  $T_0$ cuts $\Sigma$ into two pieces. We denote the piece
  corresponding to $\Delta_0$ by $\Sigma'_0$. Thus $\Sigma_0$ results
  from $\Sigma'_0$ by gluing a solid torus into its boundary, so
  $$
H_1(\Sigma_0/K_0)=H_1(\Sigma'_0/T_0)\,.$$
  Define  $\Sigma'_i$ for $i=1,\dots,n$ similarly, so
  $$
H_1(\Sigma_i/K_i)=H_1(\Sigma'_i/T_i)\,.$$
If we cut $\Sigma$ along all the tori $T_i$, $i=0,\dots,n$, the
central piece $\Sigma_{v'}$ corresponding to the node $v'$ of $\Delta$
is an $S^1$--bundle over an $(n+1)$--punctured sphere $S$. Denote a fiber
of this bundle by $f$ and the boundary components of $S$ by
$q_0,\dots,q_n$, considered as curves in $\Sigma_{v'}\cong S^1\times
S\subset \Sigma'_0$.  For $i=1,\dots,n$, $\Sigma_i$ is obtained
from $\Sigma'_i$ by gluing in a solid torus with meridian curve $f$,
so we have
  $$H_1(\Sigma_i)=H_1(\Sigma'_i)/(f)\,,$$
  so $H_1(\Sigma'_i)/(f)$ has
  order $d_i$. It follows that
  $$H_1(\Sigma'_0)/(f)=\bigoplus_{i=1}^nH_1(\Sigma'_i)/(f)$$
  has order  $d_1\dots d_n$.

  By definition of $\bar d_i$, the quotient $H_1(\Sigma'_i)/(f,q_i)$ has
  order $\bar d_i$, so the order of the element $q_i\in H_1(\Sigma'_i)/(f)$
  must be $d_i/\bar d_i$. The element $q_1+\dots+q_n\in
  H_1(\Sigma'_0)/(f)$ hence has order
  $\lcm({d_1}/{\bar d_1},\dots,{d_n}/{\bar d_n})$.

\def\ffrac#1#2{#1/#2}
  Now
  $$
H_1(\Sigma'_0/T_0)=H_1(\Sigma'_0)/(f,q_0)
  =H_1(\Sigma'_0)/(f,q_1+\dots+q_n)\,,$$
  so this group has order $\lvert
  H_1(\Sigma'_0)/(f)\rvert
  /\lcm(\ffrac{d_1}{\bar d_1},\dots,\ffrac{d_n}{\bar d_n})$. This equals
  \begin{gather*}
\ffrac{d_1\dots
  d_n}{\lcm(\ffrac{d_1}{\bar d_1},\dots,\ffrac{d_n}{\bar d_n})}=
\gcd\nolimits_{i=1}^n \bigl(\bar d_j\,\prod_{j\ne i}^n d_j\bigr)
  \end{gather*}
completing the proof.
\end{proof}
The ideal generator is also defined if $v$ is a leaf, and the above
proof shows that it equals $\vert H_1(\Sigma/K;\Z)\vert$ where $K$ is
the knot in $\Sigma$ corresponding to the leaf. Thus:
\begin{corollary}
  The order in homology of the knot in $\Sigma$ corresponding to a
  leaf of a resolution diagram $\resgraph$ is $\disc(\resgraph)/\bar d$, where
  $\bar d$ is the ideal generator at the corresponding leaf of the splice
  diagram.\qed
\end{corollary}

\section{Appendix 2: Okuma's Theorem}\label{sec:okuma}

We conjectured (Conjecture \ref{rational}) that rational singularities
and $\Q$HS-link minimally elliptic singularities are splice quotients.
In a recent preprint \cite{okuma} T.  Okuma announces, in effect, that
this conjecture is correct.  A key is an explicit construction of the
UAC, \`a la Esnault-Viehweg, via a sheaf of algebras on the resolution
of $(Y,o)$ \cite{okuma1}. The preprint \cite{okuma} is hesitant about
whether the complete intersections he constructs there, which he calls
Neumann-Wahl systems, are actually splice type. In fact, they are of
splice type. He constructs his complete intersections under a strong
condition on the graph (Condition 3.4 of \cite{okuma}).  The key point
we make is that a weaker condition, `Condition 3.3,' that he
shows this implies, is equivalent to the semigroup and congruence
conditions.  To clarify the situation we first assume this equivalence
and give versions of Okuma's main results in our language.

Recall that an \emph{end-curve} on a resolution is a rational curve
that has just one intersection point with the rest of the exceptional
divisor, so it corresponds to a leaf of the resolution graph.

\begin{theorem}[Okuma, \cite{okuma}]\label{th:end-curve}
  Let $(Y,o)$ be a normal surface singularity with $\Q$HS link whose
  resolution graph $\Gamma$ satisfies Okuma's `Condition 3.4',
  %the semigroup and congruence conditions,
  and $\bar Y\to Y$ its minimal good resolution. Suppose
  that for each end curve $E_i$ on $\bar Y$ there exists a function
  $y_i\colon Y \to \C$ such that the proper transform on $\bar Y$ of
  its zero-locus consists of one smooth irreducible curve $C_i$, which
  intersects $E_i$ transversally at one point and intersects no other
  exceptional curve.  Then $(Y,o)$ is a splice-quotient.
\end{theorem}

We describe `Condition 3.4' later; for now it suffices that it is
stronger than `Condition 3.3' (ie, the semigroup and congruence
conditions). Since it is rarely satisfied for splice-quotients, one would
prefer to replace it in the Theorem by the semigroup and congruence
conditions.  We have even conjectured that the existence of functions
$y_i$ as above is by itself equivalent to $(Y,o)$ being a
splice-quotient.  This is proved in \cite{neumann-wahl-zsphere}
(Theorem 4.1) when the link is a $\Z$--homology sphere; in this case,
the semigroup condition is deduced directly (and no UAC need be
constructed).  But this conjecture is still open in general.
\begin{proof}[Sketch of Proof of Theorem \ref{th:end-curve}]
  Let $d_i$ be the order of vanishing of $y_i$ on $C_i$. Then the
  ``Riemann surface'' of $y_i^{1/d_i}$ (ie, adjoin the $d_i$-th root
  of $y_i$ and normalize) is an abelian cover that is
  unramified away from the singular point. Thus, a $d_i$-th root $z_i$
  of $y_i$ is well defined on the universal abelian cover $(X,o)$ of
  $(Y,o)$.
If the leaves of $\Gamma$ are numbered $i=1,\dots,t$, we
  want to show that the $z_i$, $i=1,\dots,t$, embed $X$ in $\C^t$ as a
  complete intersection of splice type and $Y$ is a splice-quotient of
  $X$.

  Let us verify that the discriminant group $H_1(\Sigma)$ acts as
  specified in Proposition \ref{prop:5.3} on $\C^t$. If $\pi\colon
  \bar Y\to Y$ is the resolution, then, since the zero-divisor of
  $y_i\circ\pi$ has zero intersection with each $E_j$, this
  divisor is
$$d_iC_i -
  d_i\sum_j \bar a_{ij}E_j
%=d_iC_i + \frac{d_i}{\disc(\Gamma)}\sum_j \ell_{ij}E_j
\,,$$
where
$(\bar a_{ij})$ is the inverse of the intersection matrix
$(a_{ij})=(E_i\cdot E_j)$.
In other words,
  the order of vanishing of $y_i$ along $E_j$ is
  $-d_i\bar a_{ij}$. Thus the $d_i$-th root $z_i$
  of $y_i$ changes by $\exp(-2\pi i\bar a_{ij})$ as we go
  around a meridian curve of $E_j$. But this meridian curve represents
  the element of the discriminant group corresponding to $-e_j$ in the
  notation of Section \ref{sec:discriminant} (the sign results from
  the convention for how a fundamental group acts as covering
  transformations) so its action on $z_i$ is indeed as in Proposition
  \ref{prop:5.3} (see Lemma \ref{le:disc is inverse}).

  We next need to know that the $z_i$ generate the maximal ideal at
  our singular point. This is a significant step in Okuma's proof,
  and is where he needs `Condition 3.4' (as opposed to
  simply the semigroup and congruence conditions);
  we do not attempt to reprove it here.

  Since the semigroup and congruence conditions are satisfied by
  assumption, we can choose a system of admissible monomials in the
  $z_i$ that transform correctly under the discriminant group. Let
  $M_1,\dots,M_\delta$ be the monomials corresponding to the $\delta$
  edges at a node $v$ of the splice diagram. Then, as in the proof of
  Theorem \ref{th:change monomial} (see also the proof of Theorem 4.1
  of \cite{neumann-wahl-zsphere}), the ratios $M_i/M_\delta$ are
  invariant under the discriminant group, hence defined on $\bar Y$,
  and each has just a single zero and a single pole on $E_{v}$ (at the
  intersections with the neighboring curves $E_{v_i}$ and
  $E_{v_\delta}$ respectively). It follows that there are $\delta-2$
  linear relations among the $M_i/M_\delta$ up to higher order at
  $E_v$. Multiplying by $M_\delta$, we see that the $z_i$ satisfy a
  system of splice type equations, compatible with the discriminant
  group action. We thus get a map of $(Y,o)$ to the corresponding
  splice quotient and it is not hard to see it is an isomorphism.
\end{proof}

%Okuma proves:
\begin{theorem}[Okuma, \cite{okuma}]\label{th:okuma rational}
  If $(Y,o)$ is rational
or $\Q$HS-link minimally elliptic, then the semigroup and
congruence conditions are satisfied and functions $y_i$ as in the
above theorem exist. In addition, $(Y,o)$ is a splice quotient.
\end{theorem}

\begin{proof}
  The existence of functions $y_i$ as in Theorem \ref{th:end-curve} is
  standard in the rational case, and, as Okuma points out, follows
  also in the minimally elliptic case by the arguments of Miles Reid
  in \cite{reid} (Lemma, p.\ 112).  It thus remains to discuss
  `Condition 3.4' and the
  semigroup and congruence conditions.
%, which we do in detail, since Okuma formulates these conditions very differently. 
We first give
  some of Okuma's terminology.

A \emph{$\Q$--cycle} is a
rational linear combination of the exceptional curves $E_i$.  For each
$i$ Okuma denotes by $\bar E_i$ the $\Q$--cycle ``dual'' to $E_i$ in
the sense that $\bar E_i\cdot E_j=-\delta_{ij}$ (so $\bar E_i=-e_i$
in the notation of Section \ref{sec:discriminant}). By Theorem
\ref{th:props}:
$$\bar E_i=\frac1{\det(\Gamma)}\sum_j \ell_{ij}E_j$$
A \emph{monomial cycle} is a non-negative integer linear combination
$$D=\sum_{k\in\mathcal E} \alpha_k\bar E_k\,,\quad\text{where
}\mathcal E = (\text{the ends of }\Gamma)\,.$$ Okuma calls each
connected component $C$ of $E-E_i$, for any $i$, a \emph{branch} of
$E_i$.  We denote by $\Gamma_C$ the corresponding subgraph of
$\Gamma$.

The following is `Condition 3.4', which is well known for rational and
$\Q$HS-link minimally elliptic singularities; as Okuma says, it
follows from basic results on computation sequences, eg,
\cite{lauferrat, laufermin}).  As mentioned above, this condition is
at present needed to show the UAC is a complete intersection.

\begin{condition3.4}\label{cond:3.4}
  For any branch $C$ of any $E_i$ not an end-curve, the fundamental
  cycle $Z_C$ for $\Gamma_C$ satisfies $Z_C\cdot E_i=1$.
(The fundamental cycle is the minimal effective cycle that has
non-positive intersection number with each $E_j$ in $C$.)
\end{condition3.4}

Okuma shows this condition implies the following `Condition 3.3',
which we will show is equivalent to the semigroup and congruence
conditions:
\begin{condition3.3}\label{cond:3.3} For any node $i$ of $\Gamma$
and branch $C$ of
  $E_i$ there exists a monomial cycle $D$ such that $D-\bar E_i$ is an
  effective integral cycle (ie, with non-negative integral
  coefficients) supported on $C$. Moreover $D$ has the form
  $D=\sum_k\alpha_k\bar E_k$ with $k$ running only through the leaves
  of $\Gamma$ in $\Gamma_C$.
\end{condition3.3}
\noindent We have included the second sentence of the condition for
convenience; one can show it follows from the
first.

Okuma's proof that that `Condition 3.4' implies `Condition
3.3' is elegant and simple:
We want to add an effective integral cycle to $\bar E_i$ to get a
monomial cycle $D$ as in `Condition 3.3'.
Let $C_1=C$ be the branch in question and put $D_1=\bar E_i+Z_C$.
Clearly $D_1\cdot E_j=0$ for each $j$ outside $\Gamma_C$ and $D_1\cdot
E_j\le 0$ otherwise. If $D_1\cdot E_j=\alpha<0$ for some $j$ other than a
leaf of $\Gamma$ in $\Gamma_C$, choose such a $j$ as close as possible
to $i$, let $C_2$ be a branch of $j$ that is
in $C$, and put $D_2=D_1-\alpha Z_{C_2}$. Repeat until you have $D$ with
$D\cdot E_j=0$ for all $j$ that are not leaves of $\Gamma$ in
$\Gamma_C$. Then $D=\sum_k\alpha_k\bar E_k$ with $\alpha_k=-D\cdot
E_k$ for each $k$ a leaf of $\Gamma$ in $\Gamma_C$.
\end{proof}
%It remains to show that `Condition 3.3' is equivalent to the
%semigroup and congruence conditions.
\begin{proposition}
  `Condition 3.3' is equivalent to the semigroup and congruence conditions.
\end{proposition}
\begin{proof}
  Suppose `Condition 3.3' holds for a node $i$ of $\Gamma$. In the
  following $k$ runs through the leaves of $\Gamma$ in $\Gamma_C$ and
  $j$ runs through all vertices of $\Gamma$.
  \begin{align*}
\det(\Gamma)(D-\bar
E_i)&=\sum_k\alpha_k\sum_j\ell_{kj}E_j-\sum_j\ell_{ij}E_j\\
 &= \sum_j(\sum_k\alpha_k\ell_{kj}~-\ell_{ij})E_j\,.
  \end{align*}
  Since this cycle is supported in $C$, the coefficient of $E_i$ is
  zero: $$\sum_k\alpha_k\ell_{ki}=\ell_{ii}\,.$$ This is the semigroup
  condition.  Note that the vanishing of the coefficient of any $E_j$
  with $E_j$ not in $C$ gives a multiple of this equation, so the
  semigroup condition is equivalent to these coefficients
  vanishing for all $E_j$ not in $C$.

Now look at the coefficient of an $E_j$ that is
in $C$. The condition that $D-\bar E_i$ is an integral cycle says
$$\sum_k\alpha_k\ell_{kj}~\equiv\ell_{ij}\quad(\text{mod
}\det(\Gamma))\,.$$ Comparing with Propositions \ref{prop:congruence}
and \ref{prop:ell} we see that as $j$ runs through leaves of $\Gamma$
this is the congruence condition. Recall that the congruence condition
is an equivariance condition and the above confirms this equivariance
for the generators of the discriminant group corresponding to leaves
of $\Gamma$. The above congruence for other $E_j$ in $C$ is the
equivariance condition for the group elements corresponding to these
vertices, and thus follows once one knows it for the generators. Thus
the congruence condition for the given node and branch is equivalent to
the above congruence as $E_j$ runs through exceptional curves in $C$.
\end{proof}

\end{document}

%% file: gtspec.tex
\def\ifplaintex{\expandafter\ifx\csname documentclass\endcsname\relax}

\def\ifplaintex{\expandafter\ifx\csname documentclass\endcsname\relax}

\ifplaintex 
\else
\fi

\expandafter\ifx\csname epsfbox\endcsname\relax\input epsf\fi

\def\gt{{\mathsurround=0pt\it $\cal G\mskip-2mu$eometry \&\ 
$\cal T\!\!$opology}}        %  journal title in recommended style

\def\gtp{{\mathsurround=0pt\it $\cal G\mskip-2mu$eometry \&\ 
$\cal T\!\!$opology $\cal P\!$ublications}}  % GT publications

\def\lognumber#1{\def\thelognumber{#1}}
\def\volumenumber#1{\def\thevolumenumber{#1}}
\def\papernumber#1{\def\thepapernumber{#1}}
\def\volumeyear#1{\def\thevolumeyear{#1}}

\def\pagenumbers#1#2{\def\startpage{#1}\def\finishpage{#2}}
\def\published#1{\def\publishdate{#1}}
\def\proposed#1{\def\theproposer{#1}}
\def\seconded#1{\def\theseconders{#1}}
\def\received#1{\def\receiveddate{#1}}
\def\revised#1{\def\reviseddate{#1}}
\def\accepted#1{\def\accepteddate{#1}}

\def\asciiaddress#1{\def\theasciiaddress{#1}}
\def\asciiemail#1{\def\theasciiemail{#1}}

\long\def\asciiabstract#1{\long\def\theasciiabstract{#1}}
\def\asciikeywords#1{\def\theasciikeywords{#1}}

\def\shorttitle#1{\def\theshorttitle{#1}}

\let\\\par\let\thelognumber\relax
\let\thevolumenumber\relax\let\thepapernumber\relax
\let\thevolumeyear\relax\let\thesamplenumber\relax\let\startpage\relax
\let\finishpage\relax\let\publishdate\relax\let\receiveddate\relax
\let\reviseddate\relax\let\accepteddate\relax\let\theasciititle\relax
\let\theasciiauthors\relax\let\theasciiaddress\relax
\let\theasciiabstract\relax\let\theasciikeywords\relax
\let\theasciiemail\relax\let\theshortauthors\relax\let\theshorttitle\relax

\long\def\maketitlep{   % start of definition of \maketitlep

\count0=\startpage

\gt\hfill      %   Journal title (top left) 
\hbox to 77pt{\vbox to 0pt{\vglue -15pt\epsfbox{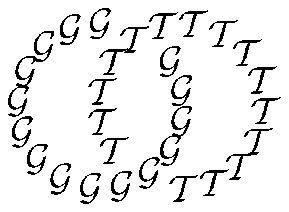}\vss}\hss}
\break
{\small\ifx\thesamplenumber\relax % sample?  
Volume \else Sample
\fi\thevolumenumber\ (\thevolumeyear)
\startpage--\finishpage\nl
Published: \publishdate}
\vglue 0.5truein plus 0.4fil minus 0.1truein

{\parskip=0pt\leftskip 0pt plus 1fil\def\\{\par\smallskip}{\ifplaintex\large
\else\Large\fi\bf\thetitle}\par\medskip}   

\vglue 0pt plus 0.1fil 

{\parskip=0pt\leftskip 0pt plus 1fil\def\\{\par}{\sc\theauthors}
\par\medskip}

\vglue 0pt plus 0.1fil 

{\small\parskip=0pt\let\newline\\
{\leftskip 0pt plus 1fil\def\\{\par}{\sl\theaddress}\par}
\expandafter\ifx\theemail\relax    % email address?
\relax\else\vglue 5pt plus 0.02fil minus 2pt\def\\{\stdspace{\rm 
and}\stdspace} 
\cl{Email:\stdspace\tt\theemail}\fi
\ifx\theurl\relax                  % URL given?
\relax\else\vglue 5pt plus 0.02fil minus 2pt\def\\{\stdspace{\rm 
and}\stdspace}
\cl{URL:\stdspace\tt\theurl}\fi\par}

\vglue 7pt plus 0.3fil minus 3pt

{\bf Abstract}
\vglue 5pt plus 0.1fil minus 2pt

\theabstract

\vglue 7pt plus 0.3fil minus 3pt

{\bf AMS Classification numbers}\quad Primary:\quad \theprimaryclass

Secondary:\quad \thesecondaryclass

\vglue 5pt plus 0.3fil minus 2pt

{\bf Keywords:}\quad \thekeywords

\vglue 10pt plus 0.5fil minus 5pt

{\small  Proposed: \theproposer\hfill Received: \receiveddate\nl
Seconded: \theseconders\hfill 
\ifx\reviseddate\relax                         % paper revised?
Accepted: \accepteddate                        % no
\else
Revised: \reviseddate                          % yes
\fi}
\eject
}       %  end of definition of \maketitlep

\font\phead=cmsl9 scaled 950
\font=cmsl9 scaled 1050
\font\pnum=cmbx10 scaled 913
\font\lnum=cmbx10 
\font\pfoot=cmsl9 scaled 950
\font=cmsl9 scaled 1050
\ifplaintex
\headline{\vbox to 0pt{\vskip -4.5mm\line{\small\phead\ifnum
\count0=\startpage ISSN 1364-0380 (on line)
1465-3060 (printed) \hfill {\pnum\folio}\else\ifodd\count0\def\\{ }% 
\ifx\theshorttitle\relax\thetitle\else\theshorttitle\fi\hfill{\pnum\folio}
\else\def\\{ and }{\pnum\folio}\hfill\ifx\theshortauthors\relax\theauthors
\else\theshortauthors\fi\fi\fi}\vss}}
\footline{\vbox to 0pt{\vglue 0mm\line{\small\pfoot\ifnum\count0=\startpage
\copyright\ \gtp\hfill\else
\gt, Volume \thevolumenumber\ (\thevolumeyear)\hfill\fi}\vss
}}
\else
\makeatletter
\def\@oddhead{{\small\lhead\ifnum\count0=\startpage ISSN 1364-0380 (on line)
1465-3060 (printed) \hfill {\lnum\number\count0}\else\ifodd\count0
\def\\{ }\ifx\theshorttitle\relax \thetitle \else\theshorttitle\fi\hfill
{\lnum\number\count0}\else\def\\{ and }{\lnum\number\count0}
\hfill\ifx\theshortauthors\relax 
\theauthors\else\theshortauthors\fi\fi\fi}}\def\@evenhead{\@oddhead}
\def\@oddfoot{\small\lfoot\ifnum\count0=\startpage\copyright\ \gtp\hfill\else
\gt, Volume \thevolumenumber\ (\thevolumeyear)\hfill\fi}
\def\@evenfoot{\@oddfoot}
\makeatother
\fi

\newwrite\gtoutfile
\long\gdef\makeheadfile{  %%% start of definition of \makeheadfile
{\def\\{, }\def\s{ }
\immediate\openout\gtoutfile head.xxx
\immediate\write\gtoutfile{Proxy-for: \ifx\theasciiauthors\relax
\theauthors\else\theasciiauthors\fi\s<\ifx\theasciiemail\relax\theemail\else\theasciiemail\fi>}
\immediate\write\gtoutfile{\noexpand\\}
\immediate\write\gtoutfile{Authors: \ifx\theasciiauthors\relax
\theauthors\else\theasciiauthors\fi}
{\def\\{ }\immediate\write\gtoutfile{Title: \ifx\theasciititle\relax
\thetitle\else\theasciititle\fi}}
\immediate\write\gtoutfile{Subj-class: GT or SG or MG etc}
\immediate\write\gtoutfile{MSC-class: \theprimaryclass\ifx\thesecondaryclass\relax\else, \thesecondaryclass\fi}
\immediate\write\gtoutfile{Journal-ref: Geom. Topol. \thevolumenumber
(\thevolumeyear) \startpage-\finishpage}
\immediate\write\gtoutfile{Comments: Published by Geometry and Topology at}
\immediate\write\gtoutfile{\s\s http://www.maths.warwick.ac.uk/gt/GTVol\thevolumenumber/paper\thepapernumber.abs.html}
\immediate\write\gtoutfile{\noexpand\\}
\immediate\write\gtoutfile{}
\ifx\theasciiabstract\relax
\immediate\write\gtoutfile{\theabstract}\else
\immediate\write\gtoutfile{\theasciiabstract}\fi
\immediate\write\gtoutfile{}
\immediate\write\gtoutfile{\noexpand\\}
\immediate\write\gtoutfile{}
\immediate\closeout\gtoutfile}}  %%% end of definition of \makeheadfile

\def\maketitlepage{\maketitlep\makeheadfile}
\let\maketitle\maketitlepage